\documentclass{article}

\usepackage{amsmath, amssymb, amsthm}
\usepackage{verbatim}
\usepackage{graphicx}

\usepackage[text={6in,9in},foot=0.6in]{geometry}

\usepackage{amsfonts}
\usepackage{graphicx} \usepackage{amsmath}

\newtheorem{theorem}{Theorem}[section]
\newtheorem{lemma}[theorem]{Lemma}

\newtheorem{proposition}[theorem]{Proposition}

\numberwithin{equation}{section} 

\newcommand{\de}{\delta}
\newcommand{\ep}{\varepsilon}

\newcommand{\ta}{\theta}
\newcommand{\mbp}{\mathbf{p}}
\newcommand{\mbq}{\mathbf{q}}
\newcommand{\mbr}{\mathbf{r}}
\newcommand{\mbx}{\mathbf{x}}
\newcommand{\mby}{\mathbf{y}}
\newcommand{\mbz}{\mathbf{z}}
\newcommand{\mbf}{\mathbf{f}}

\newcommand{\mbu}{\mathbf{u}}
\newcommand{\mbta}{\boldsymbol \theta}
\newcommand{\mbpi}{\boldsymbol \pi}

\newcommand{\enp}{e^{np_jr}}

\newcommand{\sumj}{\sum_{j=1}^n}

\newcommand{\pikb}{\pi_{kb}}

\newcommand{\e}{\varepsilon}

\begin{document}
\title{Expected Coalescence Time for a Nonuniform Allocation Process}
\author{ John K. McSweeney\footnote{Research partially supported by an NSF VIGRE fellowship. Email: {\tt mcsweeney@math.ohio-state.edu}} , \hspace{1cm} Boris G. Pittel\footnote{Research supported by an NSF grant. Email: {\tt bgp@math.ohio-state.edu}}~\footnote{Postal Address (both authors): Ohio State University, Department of Mathematics, 231 W 18th  ave, Columbus, OH 43210} \\
Department of Mathematics,  The Ohio State University}
\maketitle

\begin{abstract}
 We study a process where balls are repeatedly thrown into $n$ boxes independently
according to some probability distribution $\mbp$. We start with 
$n$ balls, and at each step all balls landing in the same box are 
fused  into a single ball;  the process terminates when there is 
only one ball left (coalescence). Let $c:=\sum_jp_j^2$, the collision probability 
of two fixed balls. We show that the 
expected coalescence time is asymptotically $2c^{-1}$, under two constraints on 
$\mbp$ that exclude a thin set of distributions $\mbp$. One of the constraints
is $c\ll \ln^{-2}n$. This $\ln^{-2} n$ is shown to be a threshold value: for $c\gg \ln^{-2} n$, there exists $\mbp$ with $c(\mbp)=c$ such that the expected coalescence time far
exceeds $c^{-1}$. Connections to coalescent processes in population biology and theoretical  computer science are discussed.\\


\end{abstract}

\medskip

\noindent{\bf Remark.~} All limits in this paper are taken as $n\to\infty$, and we use Landau notation $o,\omega, O$ in the usual sense. We  say that an event $A$ holds with high  probability (whp) if $P(A^c)=o(1)$, and we denote the set $\{1,2,\dots,n\}$ by $[n]$.

\section{Introduction} 


We consider the following  balls-into-boxes process. Let $\mbp=(p_1,\dots,p_n)$ be any 
probability vector. At time $t=0$, start with $b_0$  balls and throw them into $n$ boxes, where each ball  
has probability  $p_j$ of landing in box $j$, independently of all other balls. Fuse  all 
balls that land in the same box, and then repeat the allocation at times $t=1,2,\dots$ according  to the same rules with this possibly smaller new number of balls. The random time $T$ at which 
all balls are first fused into a single one is called the \textit{coalescence} time; we will 
be mainly interested in its asymptotic expected value.

This problem has been studied in various guises by numerous authors. Initially, 
it was stated in terms of finding the Most Recent Common Ancestor (MRCA) in a 
random genealogical process. We will restrict ourselves to the balls-into-boxes 
formulation here, but  it is important to note that any of the results we obtain 
here can be recast in the language of population biology; we will occasionally 
provide the reader the appropriate analogy. The seminal work in this area was 
done by Kingman \cite{King1, King2, King3} who, for fixed $b_0$, proved convergence of 
the underlying process to the continuous-time  \textit{coalescent} process, thereby 
establishing convergence of the distribution of $T/E[T]$.  In \cite{King1} he also 
proved that for $b_0=n$ and the uniform distribution  $\mbp=(1/n,\dots,1/n)$, 
\begin{equation}\notag
E[T] \leq 2n-2,
\end{equation}
which effectively implies that
\begin{equation}\label{Kingman2}
E[T]\sim 2n, ~~~~n \to \infty.
\end{equation}
More recently,  Donnelly and Tavar\'e \cite{DT}, M\"ohle \cite{Moehle}, and M\"ohle and Sagitov \cite{MS} have studied the limiting behavior of more general classes of allocation (reproduction) models for $b_0$ fixed. In \cite{Moehle2}, M\"ohle considered the case $b_0 \to \infty$ as well for several models, which however do not include the one at hand. 

There is a good deal of literature on continuous-time coalescent processes that bears mentioning. 
For large times $t$, when the number of balls is small, we should expect  relatively long time  intervals during  which no collisions happen, likely punctuated by binary (one-on-one) collisions. This  behavior is characteristic of Kingman-type coalescent processes, and admits a natural time scaling to a continuous-time  process. More recently, the theory of $\Lambda$- and $\Xi$-coalescents, developed by Pitman, Sagitov, and Schweinsberg \cite{Pitman, Sagitov, Schweinsberg}, among others, allows for models involving multiple simultaneous collisions. Indeed, in our model, for small times $t$, we are likely to have multiple simultaneous collisions, however these happen at {\it fixed} time intervals. There is therefore no natural time-scaling that can be performed in order to interpret this behavior in the limit as a continuous-time process with random collision times.

The process can also be described in terms of compositions of random functions: choose random functions $\{f_s:[n]\to[n]\}_{s \in \mathbf{N}}$ independently, in such a way that for all $ i\in [n]$, and for all $s\in \mathbf{N}$, $f_s(i)=j$ with probability $p_j$, independently for all $i$ and $s$. The coalescence time $T$ is then the smallest value of $t$ for which $f_t\circ\cdots \circ f_2\circ f_1$ is a constant function. It is this formulation that has been used in connection with computer science: this problem is potentially useful in bounding the running time of so-called `Coupling  from the Past' (CFTP) algorithms introduced by Propp and Wilson \cite{PW2,PW1}; see Section \ref{fut} for a brief discussion of this. Motivated by this connection, and apparently unaware of Kingman's work for the uniform distribution, Dalal and Schmutz \cite{DS} established \eqref{Kingman2};  Fill \cite{Fill} and Goh et al. \cite{GHS} derived the limiting distribution of $T/E[T]$.

For $b_0=n$, Adler et al. \cite[Theorem 4]{AAKR} were able to extend Kingman's result \eqref{Kingman2} 
to a nonuniform $\mbp$, showing that 
\begin{equation}\label{asympt}
E[T] \sim 2c_2^{-1},\quad c_2:=\sum_jp_j^2,
\end{equation}
(note that $c_2=1/n$ for the uniform $\mbp$), under the condition 
\begin{equation}\label{c3/c2}
\frac{c_3}{c_2} < \frac{3}{n},\quad c_3:=\sum_j p_j^3.
\end{equation}
Here $c_2$, $c_3$ are the probabilities of a double collision and a triple collision
respectively. In essence, \eqref{c3/c2} means that $\mbp$ is sufficiently close to
$(1/n,\dots,1/n)$. 

We should expect the largest contribution to the time $T$ to happen during the late 
stages of the process, when the number of balls is relatively small. In this case, any 
reduction in the number of balls will most likely be due to the collision  of a single pair 
of balls. This explains the appearance of $c_2$ in \eqref{asympt}.
(In fact, $c_2$ had  been used as a scaling parameter for a wide class of models 
by M\"ohle in \cite{Moehle} in the context of population genetics.) The bound \eqref{c3/c2}
on $c_3$ ensures that, if in one of those late stages the number of balls has dropped, then 
the actual decrease is exactly $1$ with conditional probability  sufficiently close to $1$.  

The proof in \cite{AAKR} revealed that the expected time spent 
in the {\it late} stages was about $2c_2^{-1}$ under conditions far 
less restrictive than \eqref{c3/c2}. \eqref{c3/c2} was used in \cite{AAKR} to show that 
the expected time  spent in the {\it early} stages was $o\big(c_2^{-1}\big)$.

In this paper we prove that this property of the process and \eqref{asympt} continue to 
hold for a much wider class of distributions $\mbp$. Here is our main result.  
\begin{theorem}\label{t1}
Let $b_0=n$. Suppose that for some $\ep >0$ however small, 
\begin{equation}\label{conds}
c_2 =o\left( \ln^{-2}n\right) ~~~\mbox{ and } ~~~~ c_3 \leq c_2^{3/2} \ln^{-(1/2+\ep)}n.
\end{equation}
Then (i)
\begin{equation}\label{result1}
E[T] = 2c_2^{-1}(1+o(1))~~~~~\mbox{for}~~n\to \infty,
\end{equation}
and (ii)
\begin{equation}\label{result2}
\frac{T}{E[T]} \Rightarrow \sum_{k \geq 2} \frac{2}{k(k-1)} Y_k,
\end{equation}
where $Y_k$ are independent and exponentially distributed, $P(Y_k >x)=e^{-x}$. 
\end{theorem} 
\noindent{\bf Remarks.}
\begin{enumerate}
\item It is always the case that $c_2 \leq 1$ and $c_3 \leq c_2^{3/2}$; so if not for the logarithmic factors, the conditions \eqref{conds} would not have excluded any probability vectors $\mbp$. Further, $c_2^{3/2}/(n^{-1}c_2)\ge n^{1/2}$, since $c_2\ge n^{-1}$ for all  $\mbp$. This means that
the restriction $c_3/c_2 < 3/n$ in \cite{AAKR} is much more stringent than 
$c_3\le c_2^{3/2}\ln^{-(1/2+
\ep)}n$ in our Theorem \ref{t1}.

\item The lower bound in \eqref{result1}, $ E[T] \geq 2c_2^{-1}(1-o(1))$,  holds for any $\mbp$ as long as $c_2\to 0$. This can be deduced from \cite[Theorem 2]{AAKR} via an elementary coupling argument; we will provide a brief proof in Section \ref{lb}. Therefore, our main task is to prove a matching upper bound.

\item For the case of the uniform $\mbp=(1/n,\dots,1/n)$, \eqref{result2} follows from 
Kingman's work. Independently, it was proved later---in a setup close to that of our paper---
by Fill \cite{Fill}; still later Goh et al. \cite{GHS} gave a detailed description of 
the cdf of the limiting distribution. 

\item Part (ii) of Theorem \ref{t1} follows from part (i), via an argument very similar to 
Fill's proof for the uniform case; cf. Theorem 6.1 of M\"ohle \cite{Moehle2}.
\end{enumerate}

Interestingly, $\ln^{-2}n$ appearing as the upper bound for $c_2$ in Theorem \ref{t1}
is a genuine threshold for the property ``$E[T]$ is of order $c_2^{-1}$ exactly'':
\begin{theorem}\label{cexthm}
Let $b_0=n$. For $c_2=\omega\left( \ln^{-2}n\right)$, there exists a probability vector $\mbp$ with 
$\sum_jp_j^2=c_2$ such that with high probability, $T =\omega\left(c_2^{-1}\right)$, and thus $E[T] =\omega\left( c_2^{-1}\right)$.
\end{theorem}
Adler et al. \cite[Theorem 5]{AAKR} had proved that $E[T]\to\infty$ for a probability vector satisfying conditions ensuring that $\lim c_2>0$. 


The rest of the paper is organized as follows. In Section \ref{earlyst}, we bound the 
expected time spent during an early phase.  We achieve this by showing that whp, for certain ``small'' values of $t$, the \textit{stochastic} process is well-approximated by  a \textit{deterministic} process, amenable to sharp estimates. To that end we will build on the method used by Pittel \cite{BGP}  for asymptotic analysis of a rumor-spreading process introduced and studied by Frieze and Grimmett \cite{FG}. In Section \ref{midlate}, we bound the expected time spent in a ``middle'' phase and a ``late'' phase, showing that the late phase contributes, overwhelmingly, to the total number of steps. In section  \ref{lb} we  prove the lower bound for Theorem \ref{t1}. In Section \ref{cexsec} we prove Theorem \ref{cexthm}. 
In the appendix we prove some auxiliary inequalities needed for the proof of Theorem \ref{t1}. 

\section{The expected duration of an early phase}\label{earlyst}

Let the distribution $\mbp$ be given. Assume that $b_0=n$. For $t \in \mathbf{N}$, we denote by $B(t)$ the random number of 
balls at time $t$, so that $B(0)=n$. In the language of the
genealogical process, $B(t)$ is the number of individuals at generation $-t$ 
which have a descendant alive in the current  generation $0$. We denote by $\tau(k)$ the random first time 
$t$ when $B(t)$ falls below $k$, i. e. $\tau(k) = \min \{t\geq 0: B(t)\leq k\}$. The 
coalescence time $T$ is therefore $\tau(1)$. Obviously, $\{B(t)\}_{t\geq 0}$ is a Markov chain 
on the state space $\{1,2,\dots,n\}$, so we will refer to $B(t)$ as the {\it state} at time $t$.

Notice that, by the definition of the stochastic sequence $\{B(t)\}_{t\geq 0}$,
\begin{equation}\notag
E[B(t+1)|B(t)=k] = \sum_{j=1}^n \left[1-(1-p_j)^k\right]~,~~~~~~~~B(0)=n,
\end{equation}
because the probability  of box $j$ receiving at least one ball out of $k$ allocated balls 
is $1-(1-p_j)^k$. It would seem natural to try to prove that the conditional 
distribution of $B(t+1)$ is concentrated around $E[B(t+1)|B(t)=k]$, as long as $k$ is
large enough. Curiously, we will be able to show instead that for smallish $t$, 
 whp $B(t+1)$ is relatively close to $\Phi_\mbp(B(t))$, where
\begin{equation}\label{phidef}
\Phi_\mbp(k) := \sum_{j=1}^n (1-e^{-p_jk}). 
\end{equation}
Note that $\Phi_\mbp(B(t))$ is close to $E[B(t+1)|B(t)]$ when most of $p_jB(t)$ are small, 
which may not be the case when $B(t)$ is relatively close to $n$.

Here is an outline of our argument. 
We introduce $k_* = o(c_2^{-1})$, and a recurrence inequality which 
the random sequence $B(t)$ is believed to satisfy whp as long as $B(t)$ is above $k_*$.
Assuming that the inequality does hold, we derive efficient bounds for $\tau(k_*)$. 
Lastly we show that indeed  whp, $\{B(t)\}_{t \leq \tau(k_*)}$ satisfies the recurrence 
inequality.    

\subsection{Variational problems}
In order to determine the likely decline of $\{B(t)\}$, we need to bound $\Phi_\mbp(k)$ from
above for a certain range of $k$. This task seems quite hard, since $\Phi_\mbp(k)$ depends
on all $n$ components of $\mbp$ in a rather complicated way. Remarkably, the worst case bounds
will do the job quite efficiently, because the worst distribution $\mbp$ turns out to
be much simpler than a feasible generic $\mbp$. 

Let us define $D_n$ to be the set of probability $n$-vectors.  For any $\mbq \in  D_n$, define 
\begin{equation} \label{Fdef}
F_\mbq(k) := \sum_{j=1}^n e^{-kq_j}, 
\end{equation}
so that  $\Phi_\mbq(k) = n-F_\mbq(k)$, and set 
\begin{equation}\notag
D(c_2):=\bigg\{\mbq \in D_n\bigg|\sum_j q_j^2 = c_2\bigg\}~~~~~\mbox{and}~~~~~D(c_2,c_3) :=\bigg\{\mbq \in D_n\bigg| \sum_j q_j^2=c_2,~\sum_j q_j^3=c_3\bigg\}.
\end{equation}
That is, $D(c_2)$ (resp. $D(c_2,c_3)$) is the set of probability vectors that share the same sum of squares (resp., squares and cubes) as $\mbp$. When dealing with $D(c_2)$, $D(c_2,c_3)$, we will assume that $c_2>1/n$; otherwise
these sets are reduced to a point $(1/n,...,1/n)$. The functional $F_\mbq(k)$ is continuous and the sets $D_n$, $D(c_2)$ and $D(c_2,c_3)$ are compact, so the infima of $F$ (as a function of $\mbq$, for fixed $k$) over these sets are attained. 
\begin{proposition}\label{var1}
For any $k \in [n]$,
\[
\min_{\mbq \in D_n} F_\mbq(k) = F_\mbu(k),
\]
where $\mbu$ is the uniform vector $\mbu=(1/n,\dots,1/n)$. 
\end{proposition}
\begin{proposition}\label{var2}
For any $k \in [n]$,
\[
\min_{\mbq \in D(c_2)} F_\mbq(k) = F_{{\mbta}(c_2)}(k),
\]
where $\mbta(c_2) = (\ta_1,\ta_2,\dots,\ta_n)$ has the property that $\ta_1>\ta_2=\dots=\ta_n$  (when listed in nonincreasing order). That is, $\mbta$ has only two distinct entries, and the larger one has support size equal to 1. We will refer to such vectors as being of ``topheavy'' type. Using the equations $\ta_1+(n-1)\ta_2=1$ and $\ta_1^2 +(n-1)\ta_2^2=c_2$, we can  explicitly express the two entries  of $\mbta$ as 
\begin{equation}\label{thetaproperties}
\theta_1=\frac{1+\sqrt{(n-1)(c_2n-1)}}{n}, ~~~~\ta_2 = \frac{1-\ta_1}{n-1}.
\end{equation}
\end{proposition}
\begin{proposition}\label{var3}
For any $k \in [n]$,
\[
\min_{\mbq \in D(c_2,c_3)} F_\mbq(k) = F_{\mbr(c_2,c_3)}(k),
\]
where $\mbr(c_2,c_3):=(r_1,r_2,\dots,r_n)$ (when ordered in nonincreasing order)  has the following  property: for  some $\nu \in [n]$, 
\begin{equation}\label{rproperties}
r_1 =\dots =r_\nu \geq r_{\nu+1} > r_{\nu+2}=\dots=r_n.
\end{equation}
That is, $\mbr(c_2,c_3)$ has  at most $3$ distinct entries, and the middle one (if any) has support size  equal to $1$.  
\end{proposition}
 Clearly, $D_n \supset D(c_2)\supset D(c_2,c_3)$, and so $F_\mbp(k)\geq F_{\mbr(c_2,c_3)}(k)\geq F_{\mbta(c_2)}(k) \geq F_\mbu(k)$. These propositions therefore provide sharper and sharper estimates, so at various junctures we will use whichever one is easiest to work with, while still being sharp enough.

\begin{proof}[Proof of Proposition \ref{var1}]
Using the fact that $\varphi(x)=e^{-x}$ is concave up,
\begin{equation}\notag
F_\mbp(k) = \sum_j e^{-p_j k}=n\sum_j \frac{1}{n}e^{-p_jk}
\geq n\exp\bigg(-\sum_j\frac{1}{n}p_jk\bigg)=ne^{-k/n}=F_\mbu(k).
\end{equation}
\end{proof}

\begin{proof}[Proofs of Propositions \ref{var2} and \ref{var3}]
To prove Proposition \ref{var2}, there are two steps:
\begin{enumerate}
\item Show that a minimizer $\mbta$ of $F_{\mbq}(k)$ on $D(c_2)$ cannot have a configuration $\ta_{j_1}>\ta_{j_2}>\ta_{j_3}$; that is, that it cannot have three distinct entries.
\item Show that a minimizer $\mbta$ cannot have a $\ta_{j_1}=\ta_{j_2}>\ta_{j_3}$ configuration ($j_i \in [n]$), which will imply that the larger entry is unique.
\end{enumerate}
We will only prove Proposition \ref{var3} as it is more difficult; the interested reader would not find it difficult to adapt the argument to prove Proposition \ref{var2}. 

We may recast Proposition \ref{var3} as follows: letting $z_j:=q_jk$, we want to minimize $G(\mbz):=\sum_j e^{-z_j}$ under the constraints 
\begin{equation}\label{e101}
\sum_j z_j=k,~~~~~ \sum_j z_j^2=k^2c_2, ~~~~~ \mbox{and}~~~~~ \sum_j z_j^3=k^3c_3. 
\end{equation}
Our task is to show that the minimizer $\mbx=(x_1,\dots,x_n)$ of $G(\mbz)$, with components listed in nonincreasing order, must have the form $x_1=\dots=x_\nu \geq x_{\nu+1} > x_{\nu+2}=\dots=x_n$ for some $\nu \in [n]$. To this end, we show first that a minimizer of $G(\mbz)$ has at most four distinct components  (Case I),  and second, that an entry value which is strictly intermediate is encountered exactly once (Case II). Our proof does not rely on the method  of Lagrange multipliers, because its applicability for the equality constraints needs a prior justification and because, in principle,  it may deliver only a ``first-order'' necessary  condition, definitely too crude to handle Case II.

{\bf Case I.~} We first show that a minimizer of $G(\mbz)$ cannot have four distinct entries. Suppose for the sake of contradiction that we have a minimizing vector $\mbx$ for which there exist $j_1,j_2,j_3,j_4$ (relabel as 1,2,3,4) such that $x_1>x_2>x_3>x_4\geq 0$. Let $y_j=x_j+\ep_j$ for $j=1,\dots,4$, and $y_j=x_j$ for $j=5,\dots,n$; we will show that for a suitable choice of $(\ep_1,\dots,\ep_4)$, $\mby=(y_1,\dots, y_n)$ satisfies the conditions \eqref{e101}, and $G(\mby)<G(\mbx)$, and thus such an $\mbx$ cannot be a minimizer on the set $D(c_2,c_3)$.

First note that we require $\ep_4 \geq 0$ because of the possibility that $x_4=0$, but $\ep_1,\ep_2,\ep_3$ can be of either sign. For $\mby$ to satisfy the conditions \eqref{e101}, we require
\begin{align}\label{ee1}
\sum_{j=1}^4 \ep_j=0,\\\label{ee2}
2\sum_{j=1}^4 x_j\ep_j+\sum_{j=1}^4 \ep_j^2=0,\\ \label{ee3}
3\sum_{j=1}^4 x_j^2\ep_j +3\sum_{j=1}^4 x_j\ep_j^2+\sum_{j=1}^4 \ep_j^3=0. 
\end{align}
Now we want $G(\mathbf{x})-G(\mathbf{y})>0$; by linearizing the $e^{-\ep_j}$ factors, it will be sufficient (by taking the $\ep_j$ as small as we wish) to show that
\begin{equation}\label{epsilons}
 e^{-x_1}\ep_1 + e^{-x_2}\ep_2+e^{-x_3}\ep_3+ e^{-x_4}\ep_4 >0.
\end{equation}
We now obtain expressions for the $\ep_j$. For given $x_j$, the system \eqref{ee1}-\eqref{ee3} is a system of $3$ \textit{nonlinear} equations in $4$ unknowns $\ep_1,\ep_2,\ep_3,\ep_4$; treating $\ep_4$ as a parameter, we hope to be able to solve it uniquely for $\ep_1,\ep_2,\ep_3$ near $(0,0,0)^T$. Let ${\boldsymbol \ep}:=(\ep_1,\ep_2,\ep_3)^T$, and write \eqref{ee1}-\eqref{ee3} as the vector equation
$$
\mbf(\boldsymbol \ep)=\bold b(\ep_4),\quad \bold b(\e_4):=(-\ep_4,-2x_4\ep_4 
-\ep_4^2,-3x_4^2\ep_4-3x_4\ep_4^2 -\ep_4^3)^T.
$$
The derivative (Jacobian) matrix of $\bold f$ at $\bold 0$ is
$$
L:=\left(%
\begin{array}{ccc}
  1 & 1 & 1 \\
  2x_1 & 2x_2 & 2x_3 \\
  3x_1^2 & 3x_2^2 & 3x_3^2 \\
\end{array}%
\right).
$$
Its determinant is equal to $6\Delta(x_1,x_2,x_3)$, where $\Delta(x_1,x_2,x_3)$ 
is the
Vandermonde determinant for $x_1,x_2,x_3$,
$$
\Delta(x_1,x_2,x_3)=(x_2-x_1)(x_3-x_1)(x_3-x_2),
$$
which is non-zero (negative), as the $x_i$ are distinct (decreasing). Therefore, by the Inverse
Vector Function Theorem (IVFT), for $|\ep_4|$ sufficiently small there exists a 
differentiable solution $\boldsymbol\ep=\boldsymbol\ep(\ep_4)$,
$\boldsymbol\e(0)=\bold 0$, such that
\begin{align*}
\boldsymbol\ep = \boldsymbol\gamma\e_4 +\bold O(\ep_4^2),~~~~
\boldsymbol\gamma:=L^{-1}\bold b^\prime(0)=L^{-1}(-1,-2x_4,-3x_4^2)^T.
\end{align*}
Explicitly, by Cramer's rule,
$$
\gamma_1=-\frac{\Delta(x_4,x_2,x_3)}{\Delta(x_1,x_2,x_3)},\quad
\gamma_2=-\frac{\Delta(x_1,x_4,x_3)}{\Delta(x_1,x_2,x_3)},\quad
\gamma_3=-\frac{\Delta(x_1,x_2,x_4)}{\Delta(x_1,x_2,x_3)}.
$$
With these formulas, \eqref{epsilons} is equivalent to showing (by letting $\ep_4>0$ be as small as needed)
$$
e^{-x_4}\Delta(x_1,x_2,x_3)<e^{-x_1}\Delta(x_4,x_2,x_3)+e^{-x_2}
\Delta(x_1,x_4,x_3)+e^{-x_3}\Delta(x_1,x_2,x_4),
$$
which in turn is equivalent to showing that
\[
D(\mbx):=
\left|%
\begin{array}{cccc}
  e^{-x_1} & e^{-x_2} & e^{-x_3} & e^{-x_4} \\
  1 & 1 & 1 & 1 \\
  x_1 & x_2 & x_3 & x_4 \\
  x_1^2 & x_2^2 & x_3^2 & x_4^2 \\
\end{array}%
\right|
>0.
\]
By using the operations typical for computation of the Vandermonde-type determinants, we get
\begin{align*}
D(\bold x)&=-\left|%
\begin{array}{ccc}
  e^{-x_2}-e^{-x_1} & e^{-x_3}-e^{-x_1} & e^{-x_4}-e^{-x_1}\\
     x_2-x_1 & x_3-x_1 & x_4-x_1 \\
   x_2^2-x_1^2 & x_3^2-x_1^2 & x_4^2-x_1^2 \\
\end{array}%
\right|\\
&=
e^{-x_1}\prod_{i=2}^4(x_i-x_1)\times
\left|%
\begin{array}{cc}
   \lambda(x_1-x_2)-\lambda(x_1-x_3) & \lambda(x_1-x_2)-\lambda(x_1-x_4)\\
      x_3-x_2 & x_4-x_2 \\
    \end{array}%
\right|,
\end{align*}
where  $\lambda(x):=(e^x-1)/x$. Next, we factor $x_3-x_2$ and $x_4-x_2$ from
the first column and from the second column. So, introducing 
\[
C(\bold x):=e^{-x_1}\prod_{i=2}^4(x_i-x_1)\prod_{i=3}^4(x_i-x_2)<0,
\]
we then get
\begin{align}\notag
D(\bold x)&=C(\bold x) \cdot
\left|%
\begin{array}{cc}
   \frac{\lambda(x_1-x_2)-\lambda(x_1-x_3)}{x_3-x_2} & \frac{\lambda(x_1-x_2)-\lambda(x_1-x_4)}{x_4-x_2}\\
      1 & 1 \\
    \end{array}%
\right|\\
&=C(\bold x) \cdot  \left[ \frac{\lambda(x_1-x_2)-\lambda(x_1-x_3)}{(x_1-x_2)-(x_1-x_3)}  - \frac{\lambda(x_1-x_2)-\lambda(x_1-x_4)}{(x_1-x_2)-(x_1-x_4)}
\right]. \label{determ}
\end{align}
Now $\lambda(x)$ is concave up for $x>0$. Therefore, since 
\[
0<x_1-x_2<x_1-x_3<x_1-x_4,
\]
the quantity in  square  brackets in \eqref{determ} is strictly negative, by considering  its  terms to be slopes of secant lines to the graph  of $\lambda(x)$. Using  this and  $C(\bold x)<0$, we get the desired conclusion, i.e. $D(\bold x)>0$.

{\bf Case II.~} Now we show that a vector $\mbx$ with a configuration $x_1>x_2=x_3>x_4\geq 0$ cannot be a minimizer of $G$ either. Define $x:=x_2=x_3$. Now that $\mbf^\prime(\mathbf{0})$ is
singular, determination of small feasible $\ep_1,\dots,\ep_4$ such that
$G(\mby)<G(\mbx)$ is more of a challenge. The fact that the linear
terms in \eqref{ee1}-\eqref{ee3} now depend on only $\ep_1,\ep_2+\ep_3,
\ep_4$ hints that $|\ep_1|, |\ep_2+\ep_3|, |\ep_4|$
should be equally small, and that $|\ep_2|$ and $|\ep_3|$, while small, should
be much larger.

Believing in this scenario, we set
$$
\ep_1=\delta_1\ep^2,\quad \ep_2=\ep+\delta_2\ep^2,\quad \ep_3=-\ep,\quad 
\ep_4=\delta_4\ep^2,
$$
and seek the feasible $\delta_i(\ep)$ for small $\ep$. To begin with, we
again require $\delta_4 \geq 0$. The conditions \eqref{ee1}-\eqref{ee3} become
\begin{align}\label{delt1}
\delta_1+\delta_2+\delta_4=& 0,\\\label{delt2}
x_1\delta_1+x\delta_2+x_4\delta_4=& -1+\ep b_2(\e,\boldsymbol\delta),\\\label{delt3}
x_1^2\delta_1+x^2\delta_2+x_4^2\delta_4 =& -2x+\ep b_3(\e,\boldsymbol\delta),
\end{align}
where the $b_i(\ep,\boldsymbol\delta)$ are polynomials. Notice that 
$\Delta:=\Delta(x_1,x,x_4)$, the determinant of the matrix in \eqref{delt1}-\eqref{delt3}, is
nonzero. So for $|\ep|$ small enough, there exists a differentiable solution
$\boldsymbol\delta(\ep)$, such that $\boldsymbol\delta(0)$ is the solution of 
\eqref{delt1}-\eqref{delt3} with $0$, $-1$, and $-2x$ respectively on the right hand side.
By Cramer's rule,
\[
\de_1(0)= \frac{1}{\Delta}
\left|%
\begin{array}{ccc}
  0 & 1 & 1  \\
  -1 & x & x_4 \\
  -2x & x^2 &  x_4^2 \\
\end{array}%
\right|,~~~
\de_2(0)= \frac{1}{\Delta}
\left|%
\begin{array}{ccc}
  1 & 0 & 1  \\
  x_1 & -1 & x_4 \\
  x_1^2 & -2x &  x_4^2 \\
\end{array}%
\right|,~~~
\de_4(0)= \frac{1}{\Delta}
\left|%
\begin{array}{ccc}
  1 & 1 & 0  \\
  x_1 & x & -1 \\
  x_1^2 & x^2 &  -2x \\
\end{array}%
\right|,
\]
which gives
\[
\de_1(0) = \frac{(x-x_4)^2}{\Delta},~~~~ \de_2(0) = \frac{(x_1-x_4)(x_1+x_4-2x)}{\Delta},~~~~ \de_4(0) = -\frac{(x-x_1)^2}{\Delta}.
\]
Reassuringly, $\de_4(0)$  is positive (because  $\Delta$ is negative). Again we want  $G(\mathbf{x})-G(\mathbf{x}+{\boldsymbol \ep})>0$; we have
\begin{align}\notag
G(\mbx)-G(\mbx+\bold\ep) =&\sum_{j=1}^4 (1-e^{-\ep_j})e^{-x_j}\\ \notag
& e^{-x_1}(1-e^{-\ep_1}) +e^{-x}(1-e^{-\ep_2}+1-e^{-\ep_3})+e^{-x_4}(1-e^{-\ep_4})\\ \notag
=&e^{-x_1}(\de_1\ep^2+O(\ep^4)) \\ \notag
&+e^{-x}(\ep+\de_2\ep^2 -(1/2)(\ep+\de_2\ep^2)^2+O(\ep^3)-\ep -(1/2)\ep^2 +O(\ep^3))\\ \notag
& +e^{-x_4}(\de_4\ep^2+O(\ep^4))\\ \label{deltazero}
=&\ep^2 (\de_1e^{-x_1}+(\de_2-1)e^{-x}+\de_4e^{-x_4}+O(\ep)).
\end{align}
Thus, by taking $\ep$ sufficiently small, \eqref{deltazero} will  be $>0$  if
\[
\de_1(0) e^{-x_1} +(\de_2(0)-1)e^{-x} + e^{-x_4} \de_4(0) >0.
\]
In light of the formulas for the $\delta_i(0)$ and the fact that $\Delta<0$, this is equivalent to
\begin{equation}\label{tx}
T(\mbx) := -(x-x_4)^2 e^{-x_1} +(\Delta - (x_1-x_4) (x_1+x_4-2x)) e^{-x}+ (x-x_1)^2 e^{-x_4} >0.
\end{equation}
Now multiplying $T(\mbx)$ by $e^x$ and using the inequalities 
$$
e^{-(x_1-x)}< 1-(x_1-x)+(x_1-x)^2/2,\quad e^{x-x_4}> 1+ 
(x-x_4)+(x-x_4)^2/2,
$$
we get
\begin{align*}
e^xT(\mbx)&=(x-x_4)^2 (1-e^{x-x_1}) + (x-x_1)^2(e^{x-x_4}-1) +\Delta\\
&> (x-x_4)^2\left((x_1-x)-\frac{(x_1-x)^2}{2}\right) + (x-x_1)^2\left((x-x_4) +\frac{(x-x_4)^2}{2}\right) +\Delta\\
&=(x-x_4)^2(x_1-x)+(x-x_1)^2(x-x_4)+(x-x_1)(x_4-x_1)(x_4-x)\\
&= 0.\quad (!)
\end{align*}
Therefore \eqref{tx} holds, and thus as before, $\mbx$ cannot  be a minimizer, and this concludes Case II. 

This only  leaves the possibility that  the minimizer $\mbx$  of $G$ is of the  form
\[
x_1=x_2=\dots=x_\nu\geq x_{\nu+1} \geq x_{\nu+2} =\dots=x_n,
\]
and thus that the minimizer $\mbr$ of $F_\mbq(k)$ over $D(c_2,c_3)$ is of the form \eqref{rproperties}, for any  $k$. This concludes the proof of Proposition \ref{var3}.
\end{proof}

\subsection{Identifying and iterating a likely recurrence inequality}

Let 
\begin{equation}\label{kstardef}
k_* := c_2^{-1}\ln^{-\ep}n,
\end{equation}
with $\ep$ coming from \eqref{conds}. ($k_*$ is meant to be an integer, as is  another
parameter $k_1$ defined later, but for simplicity we omit the ``integer part'' notation.)
This $k_*$ will serve as a threshold separating the ``early'' states ($B(t)> k_*)$ from
the ``late'' states ($B(t)\le k_*$). So, in light of the informal discussion
in the introduction, ``$k_*=o(c_2^{-1})$'' should be more or less expected; the need 
for an additional factor, $\ln^{-\ep}n$, will become clear later, in Section 3.

Our immediate task is to identify a function $\Psi_{\mbp}(k)$, such that, intuitively
at least, the random sequence $\{B(t)\}$  whp satisfies a recurrence inequality
\begin{equation}\label{rec10}
B(t+1)\le \Psi_{\mbp}(B(t)),\quad \mbox{if }\quad B(t)\ge k_*.
\end{equation}
Then, for $k\ge k_*$, $\Psi_{\mbp}(k)$ needs to be large enough so that, conditionally
on $\{B(t)=k\}$, the event $\{B(t+1)\le \Psi_{\mbp}(B(t))\}$ is very likely. Since we think
that $\Phi_{\mbp}(k)$ defined in \eqref{phidef} is a ``sharp'' conditional predictor of $B(t+1)$, 
we must have $\Psi_{\mbp}(k)> \Phi_{\mbp}(k)$. Also, to be of any use, $\Psi_{\mbp}(k)$
must fall below $k$. Last, but not least, we must be able to solve a chosen recurrence.
The function 
\begin{equation}\label{psidef}
\Psi_{\bold p}(k)=\frac{1}{2}(k+\Phi_{\bold p}(k))=
\frac{1}{2}\left(k+n-F_{\bold p}(k)\right),
\end{equation}
with $F_\mbp(k)$ as defined in \eqref{Fdef}, happens to meet all three requirements.

Define an event
$$
\Delta:=\{\forall\, t\,:\,B(t)\ge k_*\Longrightarrow B(t+1)\le \Psi_{\bold p}(B(t))\}.
$$
In Section 2.5 we will show that $P(\Delta)\to 1$ as $n\to\infty$. In this section, assuming
that the event $\Delta$ holds, we solve the recurrence \eqref{rec10} and estimate sharply
$\tau(k_*)$, the first moment $t$ when $B(t)\le k_*$.

\begin{lemma}\label{deltalemma}
On the event $\Delta$, 
\[
\tau(k_*) \le 5c_2^{-1/2}\ln n=o\big(c_2^{-1}\big).
\]
\end{lemma}

\begin{proof}The proof is divided into two cases.

{\bf Case I.~} $c_2 \geq 2n^{-1}$. First of all, by Proposition \ref{var2}, we have 
\begin{equation}\notag
\Psi_\mbp(k) \leq \Psi_{\mbta(c_2)}(k)~~~~~\forall k \in [n],
\end{equation}
$\mbta(c_2)$ being the topheavy distribution $(\ta_1,\ta_2,\dots,\ta_2)$ with parameter $c_2$. Therefore on the event $\Delta$, 
\begin{equation}\label{thetacor1}
B(t+1) \leq \Psi_{{\mbta}(c_2)}(B(t)),~~~~~\forall t \leq \tau(k_*).
\end{equation}

Let us bound $\Psi_{\mbta(c_2)}(B(t))$ from above. Since $c_2 \geq 2n^{-1}$, we get, using \eqref{thetaproperties},
\begin{equation}\label{thetaineq}
\ta_1 \geq \frac{1}{2}c_2^{1/2}~~\Rightarrow~~\ta_2 \leq \frac{1-(1/2)c_2^{1/2}}{n-1}.
\end{equation}
By \eqref{thetaineq} and 
\[
\Psi_{\mbta}(B(t))=\frac{1}{2}(B(t)+n-e^{-\ta_1 B(t)} -(n-1)e^{-\ta_2B(t)}),
\]
we have (using the inequality $e^{-x}\geq 1-x$)
\begin{align}\notag
B(t+1) &\leq \frac{1}{2}\left(B(t)+n -0-(n-1)(1-\ta_2B(t))\right)\\\label{Nonunifrec1}
& \leq \frac{1}{2} \left(B(t)+(1-c_2^{1/2}/2)B(t)+1\right)\\\label{Nonunifrec}
&\leq \left(1-\frac{c_2^{1/2}}{4}\right) B(t) +\frac{1}{2},
\end{align}
a \textit{linear} recurrence inequality. (Implicit in this derivation is an intuition that, for the distribution $\mbta$ in question, a large enough proportion of collisions happen in box 1, and that we may disregard collisions in boxes $2,\dots,n$ without inducing too large an error.) It follows from \eqref{Nonunifrec} and $B(0)=n$ that 
\begin{equation}\label{solution}
B(t) \leq n\left(1-\frac{c_2^{1/2}}{4}\right)^t +2c_2^{-1/2}~~~~~~~\mbox{for}~~B(t) \geq k_*.
\end{equation}
To get a bound on $\tau(k_*)$, notice that
\begin{equation}\label{Nonunifrec2}
k_* < B(\tau(k_*)-1) \leq n \left(1-\frac{c_2^{1/2}}{4}\right)^{\tau(k_*)-1} +2c_2^{-1/2}
\end{equation}
and let $\tau:=\tau(k_*)-1$. Now $c_2^{-1/2} =o(k_*)$ if $c_2 = o(\ln^{-2\ep}n)$, which is certainly implied by \eqref{conds}. So we can crudely use the bound $2c_2^{-1/2}\leq (1/2)k_*$ in \eqref{Nonunifrec2} to get
\begin{equation*}
(1/2)k_* \leq n(1-c_2^{1/2}/4)^\tau \leq n\exp\left(-c_2^{1/2}\tau/4\right). 
\end{equation*}
Taking logarithms and solving for $\tau$, we obtain
\[
\tau \leq 4\left(\ln n +\ln c_2 +\ep\ln\ln n +\ln 2\right)c_2^{-1/2} \leq 5c_2^{-1/2}\ln n.
\]
Hence $\tau(k_*) =o(c_2^{-1})$, since $c_2 =o(\ln^{-2}n)$, which is the first  condition in \eqref{conds}. 

{\bf Case II.~} $c_2 \leq 2n^{-1}$. This time $\mbta(c_2)$ is too close to being uniform, and the inequality \eqref{Nonunifrec1} is too crude. A bit of reflection shows that we should not expect that $B(t)$ decay exponentially here. We show instead that, for some absolute constant $A$,
\begin{equation}\label{eq1}
B(t)\leq \frac{An}{t+1}, ~~~~~t\leq \tau(k_*).
\end{equation}
The  proof is by induction. The case $t=0$  holds if $A\geq 1$. Suppose \eqref{eq1} holds  for some  $t$. Observe that 
\[
\Psi_\mbp(k) =\sum_{j=1}^n \psi(p_jk),~~~~~\psi(x):= \frac{1}{2} (x+1-e^{-x}),
\] 
and that $\psi(x)$ is increasing and concave down. Then by the inductive assumption,
\begin{align*}
B(t+1) \leq \sum_{j=1}^n \psi(p_jB(t)) \leq \sum_{j=1}^n \psi\left(p_j \frac{An}{t+1}\right)
\leq n\psi\left(\frac{1}{n}\sum_{j=1}^n p_j\frac{An}{t+1}\right) =n\psi\left(\frac{A}{t+1}\right). 
\end{align*}
So we need to find $A \geq 1$ such that 
\begin{equation}\label{induct}
\psi\left(\frac{A}{t+1}\right) \leq A/(t+2),
\end{equation}
or, defining $x:=A/(t+1)$,
\begin{equation}\label{eq1a}
 \psi(x) \leq x~\frac{t+1}{t+2}~~\Longleftrightarrow 1-e^{-x} \leq x \left(1-\frac{2}{t+2}\right).
\end{equation}
We therefore define $x(t)$ implicitly by
\begin{equation}\label{eq2}
1-e^{-x(t)} = x(t)\left(1-\frac{2}{t+2}\right);
\end{equation}
by considering the graphs of the functions of $x(t)$ on the left- and right-hand sides of \eqref{eq2}, it is clear that  $x(t)$  is well-defined and decreasing (to 0)  in $t$. Therefore \eqref{eq1a} is satisfied iff $x \geq x(t)$. It is not difficult to show that   
\[
x(t) \sim \frac{4}{t+2},~~~~t\to \infty,
\]
and so $A_*:= \limsup_{t\to\infty} (t+1)x(t)$ is finite.  Thus to satisfy \eqref{induct} and thereby to complete the inductive proof, we can pick $A=\max\{1,A_*\}$.  

Therefore, on the event $\Delta$, we have
\[
k_* < B(\tau(k_*)-1) \leq \frac{An}{\tau(k_*)}
\]
which we  can invert to get
\[
\tau(k_*) \leq \frac{An}{k_*} =2A\ln^\ep n,
\]
which is certainly $o\bigl(c_2^{-1/2}\ln n\bigr)$. Combining this with the case $c_2 \geq 2/n$, we have $\tau(k_*)=O\bigl(c_2^{-1/2}\ln n\bigr)$ on the event $\Delta$.  This completes the proof of Lemma \ref{deltalemma}.
\end{proof}


\subsection{Exponential tail bounds}
We need to show that $P(\Delta)$ converges to $1$, and that it does so sufficiently fast.

To this end, and also for Theorem \ref{cexthm}, we establish two-sided exponential tail bounds for for the  distribution  of $B(t+1)$ conditioned on $B(t)$. Like Chernoff bounds for sums of i.i.d. random variables, the bounds are based on a generating function approach for estimating the probabilities of
large deviations.  

Let $\pi_{kb}:= P(B(t+1)=b|B(t)=k)$, and note that $\pi_{kb}=0$ for $k<b$. Introduce
the tail probabilities
\begin{equation}\notag
\pi_{kb}^- = P(B(t+1)<b|B(t)=k)~~~~~\mbox{and}~~~~~~~\pi_{kb}^+ =
P(B(t+1)>b|B(t)=k).
\end{equation}

\begin{theorem}\label{Chern1}
\begin{equation}\label{Chernoff1}
\pi_{kb}^- \leq 3\sqrt{k}
\exp\left[-\frac{(\Phi_\mbp(k)-b)^2}{2k}\right],~~~~~ b\leq
\Phi_\mbp(k)
\end{equation}
and
\begin{equation}\label{Chernoff2a}
\pi_{kb}^+ \leq 3\sqrt{k}
\exp\left[-\frac{(b-\Phi_\mbp(k))^2}{2k}\right],~~~~~b\geq
\Phi_\mbp(k).
\end{equation}
\end{theorem}

\begin{proof}[Proof of Theorem \ref{Chern1}]

The heart of the proof is an expression for $\pi_{kb}$ by means of generating functions. As usual, the expression
$[y^m]f(y)$ denotes the coefficient of $y^m$ in the power series expansion of $f(y)$.

\begin{lemma}\label{gf}
\begin{equation*}
\pi_{kb}=\frac{k!}{n^k} [x^kz^b] \left(\prod_{j=1}^n
(1+z(e^{np_jx}-1))\right),~~~~~~1\leq b\leq k \leq n.
\end{equation*}
\end{lemma}
\begin{proof}

\begin{align*}
\pi_{kb} & = P(B(t+1)=b|B(t)=k)\\ &= \sum_{U\subset
    [n],~|U|=b} P(\mbox{$k$ balls go into exactly the boxes indexed
  by $U$})\\ &= \sum_U \mathop{\sum_{\ep_1+\dots+\ep_b=k}}_{
  \ep_j>0~\forall j \in [b]} P(\mbox{each box $j$ from $U$ gets
  $\ep_j$ balls})\\ &=\sum_U \sum_{\vec{\ep}} \binom{k}{\ep_1 \ep_2
  \cdots \ep_b} \prod_{j\in U} p_j^{\ep_j}\\ &=k! \sum_U
\sum_{\vec{\ep}} \prod_{j\in U} \frac{p_j^{\ep_j}}{\ep_j!}.
\end{align*} 
Now we build a {\it bivariate\/} generating function for the probabilities $\pi_{kb}$. Start with
the $k$ index. Incorporating for future convenience an $n^k$ factor, we have
\begin{align*}
\sum_k \pi_{kb} \frac{n^k}{k!}x^k &= \sum_U \sum_k
\mathop{\sum_{\ep_1+\dots+\ep_b=k}}_{\ep_j >0} \prod_{j\in
  U}\frac{(nx)^{\ep_j}p_j^{\ep_j}}{\ep_j!}.
\end{align*}
We merge the second and third sums, yielding for the right-hand
side
\[
\sum_U \sum_{\ep_j >0} \prod_{j\in
  U}\frac{(nx)^{\ep_j}p_j^{\ep_j}}{\ep_j!}.
\]
Reversing the order of summation and multiplication, we get
\begin{align*}
\sum_k \pi_{kb} \frac{n^k}{k!}x^k=&\sum_U \prod_{j\in U}
\sum_{\ep=1}^\infty \frac{(nxp_j)^{\ep}}{\ep!}\\ =&\sum_U \prod_{j\in
  U} (\exp(np_jx)-1).
\end{align*}
Multiplying by $z^b$ and summing for $b \geq 0$, we obtain
\begin{align*}
\sum_b \sum_k \pi_{kb} \frac{n^k}{k!}x^k z^b &= \sum_b z^b
\sum_{|U|=b} \prod_{j\in U} (\exp(np_jx)-1)\\ &= \sum_{U\subset [n]}
z^b \prod_{j\in U} (\exp(np_jx)-1)\\ &= \sum_{U\subset [n]}
\prod_{j\in U} z(\exp(np_jx)-1)\\ &= \prod_{j=1}^n \left[
  1+z(\exp(np_jx)-1)\right].
\end{align*}
Therefore we have
\[
\pi_{kb} \frac{n^k}{k!} = [x^kz^b] \prod_{j=1}^n \left[
  1+z(\exp(np_jx)-1)\right]
\]
and from here the lemma follows.
\end{proof}

By Lemma \ref{gf},
\[
g_k(z):=E[z^{B(t+1)}|B(t)=k]= \sum_{b=1}^k \pi_{kb}z^b=
\frac{k!}{n^k}[x^k] \prod_{j=1}^n (1+z(e^{np_jx}-1)).
\]
Now, for $b\leq k$, we have
\begin{align*}
z^b\pikb^- = \sum_{i=1}^b \pi_{ki}z^b \leq ~\sum_{i=1}^b \pi_{ki} z^i
\leq~ g_k(z),~~~~ 0<z\leq 1,\\ z^b\pikb^+ = \sum_{i=b+1}^k \pi_{ki}z^b
\leq \sum_{i=b+1}^k \pi_{ki} z^i \leq g_k(z),~~~~ z\geq 1.~~~~
\end{align*}
This gives
\begin{align}\label{firstgf1}
\pi_{kb}^- \leq \frac{g_k(z)}{z^b}=\frac{k!}{z^bn^k}[x^k]
\prod_{j=1}^n (1+z(e^{np_jx}-1)),~~~~~~\mbox{for}~~0<z \leq 1,\\
\label{firstgf2}
\pi_{kb}^+ \leq \frac{g_k(z)}{z^b}=\frac{k!}{z^bn^k}[x^k]
\prod_{j=1}^n (1+z(e^{np_jx}-1)),~~~~~~\mbox{for}~~z \geq 1.~~~~~
\end{align}
Since the coefficients of the products in
\eqref{firstgf1}-\eqref{firstgf2} are nonnegative, we use the
inequality $[x^k]f(x) \leq f(x)/x^k ~~(\forall~x>0)$ to obtain
\begin{align}\label{prod1}
\pikb^- \leq \frac{k!}{z^b (nr)^k}
\prod_{j=1}^n(1+z(e^{np_jr}-1))~,~~~\forall ~r>0,~~0<z \leq1,\\
\label{prod2}
\pikb^+ \leq \frac{k!}{z^b (nr)^k}
\prod_{j=1}^n(1+z(e^{np_jr}-1))~,~~~\forall ~r>0,~~z\geq 1.~~~~~~
\end{align}
Our task is to get the most out of these bounds
\eqref{prod1}-\eqref{prod2} by choosing values for $z$ and $r$
judiciously. We use Stirling's formula $k! \leq 3 \sqrt{k} (k/e)^k$ to
transform the product-type formulas \eqref{prod1} and \eqref{prod2}
into
\begin{align}
\pikb^- \leq 3\sqrt{k} \exp(H(z,r,b)) \label{q1bound}~,~~~~z \leq 1,\\ 
\pikb^+ \leq 3\sqrt{k} \exp(H(z,r,b)) \label{q2bound}~,~~~~z \geq 1,
\end{align}
where
\begin{equation}\label{H}
    H(z,r,b):=k\ln\left(\frac{k}{rne}\right)-b\ln(z)+\sum_j
    \ln\bigl(1+z(e^{np_jr}-1)\bigr).
\end{equation}
For a given $b$, we want to use  a stationary point  of
$H(z,r,b)$, i.e. a solution to
\begin{align}\label{hz}
H_z = -\frac{b}{z} + \sum_j \frac{e^{np_jr}-1}{1+z(e^{np_jr}-1)}=0,\\
\label{hr}
H_r =-\frac{k}{r} +z \sum_j \frac{np_j e^{np_jr}}{1+z(e^{np_jr}-1)}
=0.
\end{align}
This complicated system has a simple solution $(z_*,r_*)=(1,k/n)$ for
\[
b=b_*:= n-\sum_je^{-p_jk}\,\, (=\Phi_{\mbp}(k)).
\]
Moreover, from \eqref{H} it is immediate that $H(z_*,r_*,b_*)=0$. This
is a first sign that the inequalities \eqref{q1bound}-\eqref{q2bound}
may indeed lead to meaningful explicit bounds for $\pi^{\pm}_{kb}$. Of course,
 we need to know that \eqref{hz}-\eqref{hr} has a solution $(z,r)$  for
$b\neq b_*$ as well, such that $z< 1$ for $b<b_*$, and $z>1$ for $b>b_*$. 

Observe that
\[
\mbox{det} \left(%
\begin{array}{cc}
  H_{zz} & H_{zr} \\ H_{rz}&H_{rr} \\
\end{array}%
\right) = H_{zz}H_{rr} - H_{zr}^2 >0,
\]
for every solution $(z,r)$ of \eqref{hz}-\eqref{hr} (see Lemma
\ref{hess} in the appendix for a proof of this).  So, by the implicit
vector function theorem, there exists an infinitely differentiable
solution $(z(b),r(b))$ of \eqref{hz}-\eqref{hr}, such that
\[
(z(b_*),r(b_*))=(z_*,r_*).
\]
Moreover, $z(b)$ is strictly increasing (see \eqref{derz} in the
appendix), so that indeed $z(b)<1$ for $b<b_*$, and $z(b)>1$ for
$b>b_*$.  So, introducing
\begin{equation}\label{hdef}
h(b)=H(z(b),r(b),b),
\end{equation}
we have
\begin{align}
\pi_{kb}^{-}\le& 3\sqrt{k}\exp(h(b)),\quad
b<b_*,\label{e47}\\ \pi_{kb}^{+}\le& 3\sqrt{k}\exp(h(b)),\quad
b>b_*;\label{e48}
\end{align}
here $h(b_*)=H(z_*,r_*,b_*)=0$.

To get efficient bounds from \eqref{e47}-\eqref{e48}, let us
approximate $h(b)$ by its Taylor polynomial about $b_*$. First, using
\eqref{hz}-\eqref{hr},
\begin{align}\notag
h^\prime(b)&=\frac{d}{db}H(z(b),r(b),b)\\\notag
&=H_z(z(b),r(b),b)z^\prime(b)+
H_r(z(b),r(b),b)r^\prime(b)+H_b(z(b),r(b),b)\\\label{hprime}
&=H_b(z(b),r(b),b)=-\ln(z(b)).
\end{align}
It follows that $h(b)$ is unimodal (concave down, in fact), attaining
its zero maximum at $b=b_*$. Consequently
\begin{equation}
h(b)=h(b_*)+h^\prime(b_*)(b-b_*)+\frac{h^{\prime\prime}(\tilde
  b)}{2}\, (b-b_*)^2=\frac{h^{\prime\prime}(\tilde b)}{2}\,
(b-b_*)^2,\label{e49}
\end{equation}
$\tilde b$ being between $b$ and $b_*$. It is shown in the appendix
(Lemma \ref{hdplemma}) that
\[
h^{\prime\prime}(\tilde b)\le -\frac{1}{k}.
\]
This bound and \eqref{e47}, \eqref{e48}, \eqref{e49} imply
\eqref{Chernoff1} and \eqref{Chernoff2a}, thereby concluding the proof of Theorem \ref{Chern1}.
\end{proof}

\subsection{Using the exponential tail bounds}\label{etbds}

For the upper bound, we will need \eqref{Chernoff2a}, which gives, for
$b=\Psi_\mbp(k)$ (recall the definition \eqref{psidef} of
$\Psi_\mbp(k)$),
\begin{align}\label{e100}
\pi_{k,\Psi_\mbp(k)}^+ &\leq 3\sqrt{k}
\exp\left[-\frac{(\Psi_\mbp(k)-\Phi_\mbp(k))^2}{2k}\right]\\ &=
3\sqrt{k} \exp\left[-\frac{(k-\Phi_\mbp(k))^2}{8k}\right].
\end{align}
 Introducing
\[
H_\mbp(k) := \frac{1}{k}(\Psi_\mbp(k)-\Phi_\mbp(k))^2,
\]
we rewrite \eqref{e100} as
\begin{equation}\label{e104}
\pi_{k,\Psi_\mbp(k)}^+ \leq 3\sqrt{k}
\exp\left(-\frac{1}{2}H_\mbp(k)\right).
\end{equation}
The next Lemma states, roughly, that the larger $k$ is, the more
likely it is that the next state $B(t+1)$ is close to the prediction
based on information $B(t)=k$.
\begin{lemma}\label{decr}
For all $k$ and $\mbp$, $H_\mbp(k)$ is increasing in $k$.
\end{lemma}
\begin{proof}
Define
\begin{equation}\notag
N_\mbp(k):= \Psi_\mbp(k)-\Phi_\mbp(k) = \frac{1}{2}\left(k-n+\sumj
e^{-p_jk}\right)
\end{equation}
so that $H_\mbp(k)=k^{-1}N_\mbp(k)^2$. Then
$$
H_\mbp^\prime(k)=\frac{N_\mbp(k)}{k^2}\big[2kN_\mbp^\prime(k)-N_\mbp(k)\big],
$$ 
where
$$
2kN_\mbp^\prime(k)-N_\mbp(k)=\frac{1}{2}\sum_{j=1}^n\big(p_jk+1-(2p_jk+1)e^{-p_jk}\big)>0,
$$ 
because
$$ 
f(x):=x+1-(2x+1)e^{-x}>x+1-\frac{2x+1}{x+1}=\frac{x^2}{x+1}>0 \mbox{~~for~~} x>0.
$$
This completes the proof of Lemma  \ref{decr}.
\end{proof}

\begin{lemma}\label{Hlemma}
For all $k \in \{k_*,\dots,n\}$,
\[
H_\mbp(k) \geq A_*\ln^{1+\ep}n,
\]
where $A_*$ is some absolute constant. Thus the probability in
\eqref{e104} is \emph{superpolynomially} small.
\end{lemma}

\begin{proof}
Consider first the case $c_2\ge 2 n^{-1}$. From Proposition \ref{var3}
and Lemma \ref{decr}, it follows that, for all $k\ge k_*$, $H_{\mbp}(k) \ge H_{\mbp}(k_*) \ge H_{\mbr}(k_*)$; here
\begin{align}\notag
H_\mbr(k_*) &= (\Psi_\mbr(k_*)-\Phi_\mbr(k_*))^2/k_*\\ \notag
&=\frac{1}{4k_*}\left[\sum_{j=1}^n
  \big(e^{-r_jk_*}-1+r_jk_*\big)\right]^2 \\ \label{hmbrk} &\geq
\frac{1}{4k_*}\bigg[\nu \big(e^{-r_1k_*}-1+r_1k_*\big)\bigg]^2.
\end{align}
To bound \eqref{hmbrk} from below we need to have sharp bounds for
$r_1$ and $\nu$. Recalling the definition of $\bold r$ in Lemma
\ref{var3}, and letting $\mu:=n-\nu-1$, we have
\begin{equation}\label{e58}
\nu r_1+r_2+\mu r_3=1,\quad \nu r_1^2+r_2^2+\mu r_3^2=c_2,\quad \nu
r_1^3+r_2^3+\mu r_3^3=c_3.
\end{equation}
Obviously $r_3\le n^{-1}$. Since we assume that $c_2 \ge 2n^{-1}$, we
also have
\[
c_3 \ge c_2^2\ge \frac{4}{n^2}.
\]
Hence
\begin{equation}\label{e59}
\mu r_3^2\le \frac{\mu}{n^2}\le \frac{1}{n}\le \frac{c_2}{2},\quad \mu
r_3^3\le \frac{\mu}{n^3}\le \frac{1}{n^2}\le \frac{c_3}{4}.
\end{equation}
Combining \eqref{e58}, \eqref{e59} and $r_2 \leq r_1$ we get
\begin{equation}\label{m11}
\frac{c_2}{4}\le \nu r_1^2\le c_2,\quad \frac{3c_3}{8}\le \nu r_1^3\le
c_3.
\end{equation}
These double inequalities imply directly that
\begin{equation}\label{e60}
\frac{3}{8}\frac{c_3}{c_2}\le r_1 \le 4\frac{c_3}{c_2},\quad
\frac{1}{64}\frac{c_2^3}{c_3^2}\le \nu \le
\frac{64}{9}\frac{c_2^3}{c_3^2}.
\end{equation}

Armed with \eqref{e60} we return to \eqref{hmbrk}. Recalling that
$k_*=c_2^{-1}\ln^{-\ep}n$, we need to consider separately the
subsequences $\{n_i\}$ such that $r_1k_*=O(1)$ for $n\in \{n_i\}$, and
the subsequences $\{n_i\}$ along which $r_1k_*\to\infty$. In the first
case
\[
\inf_{n_i}\frac{e^{-r_1k_*}-1+r_1k_*}{r_1^2k_*^2}\ge \delta >0.
\]
So, using \eqref{m11}, we obtain from \eqref{hmbrk}:
\[
H_{\bold r}(k_*)\ge \frac{\delta^2}{4}\nu^2 r_1^4 k_*^3\ge
\frac{\delta^2}{4}(\nu r_1^2)^2c_2^{-3}\ln^{-3\ep}n \geq \frac{\delta^2}{64} c_2^{-1}\ln^{-3\ep}n.
\]
In the second case, for $n_i$ large enough,
\[
e^{-r_1k_*}-1+r_1k_*\ge \frac{1}{2} r_1k_*.
\]
So, using \eqref{m11} again and \eqref{e60},
\[
H_{\bold r}(k_*)\ge \frac{1}{16}\nu^2r_1^2k_*=\frac{1}{16}\nu(\nu r_1^2)k_* \geq
\frac{1}{2^{12}}\frac{c_2^3}{c_3^2}\ln^{-\ep}n.
\]
Therefore there exists a constant $A>0$ such that, for $n$ large
enough,
\[
H_{\bold r}(k_*)\ge A
\min\left\{c_2^{-1}\ln^{-3\ep}n,\,\frac{c_2^3}{c_3^2}\ln^{-\ep}
n\right\},
\]
for if there weren't, there would  be  a subsequence $\{ r_1k_*(n_j)\}_{j\geq1}$ along which this did not hold. This subsequence would then have a further subsequence which tended to infinity, or remained bounded, which would contradict  the bound for one of the two cases established above.  
Since
\begin{equation*}
c_3\le c_2^{3/2}\ln^{-(1/2+\ep)}n,\quad c_2\le C\ln^{-2}n
\end{equation*}
for some large enough constant $C$, the last inequality leads to
\[
H_{\bold r}(k_*)\ge A \ln^{1+\varepsilon}n,
\]
for $n$ large enough, as long as $\ep<1/4$, which of course  we may assume without  loss of generality.\\

Consider now the case $c_2\le 2 n^{-1}$. This time, by Proposition
\ref{var1} and Lemma \ref{decr}, for all $k\ge
k_*=c_2^{-1}\ln^{-\ep}n$,
\[
H_{\mbp}(k)\ge H_{\mbp}(k_*)\ge H_{\bold u}(k_*),
\]
where
\begin{align*}
H_{\bold u}(k_*)=\frac{1}{4k_*}\left(k_*-n+ne^{-k_*/n}\right)^2 \ge
\frac{k_*^3}{36 n^2}\ge \frac{1}{288} n\ln^{-3\ep}n.
\end{align*}
(For the first inequality we used
\[
e^{-x}-1+x \ge \frac{x^2}{2}-\frac{x^3}{6}\ge \frac{x^2}{3},\quad \mbox{for}~
x\in (0,1).)
\]
Thus $H_{\mbp}(k)=\omega\left( \ln^{1+\varepsilon} n\right)$ for this case, and this
concludes the proof of Lemma \ref{Hlemma}.
\end{proof}

\subsection{Validation of the deterministic approximation} \label{validation}

Now that we have established the (superpolynomially) small bound for
$\pi^{+}_{k,\Psi_{\mbp}(k)}$, we can finally show that the event
$\Delta$ is extremely likely indeed.

\begin{lemma}\label{lemma3.8}
 For some constant $A>0$,
\[
P(\Delta)\ge 1 - 12n\ln n\exp\left(-A\ln^{1+\varepsilon} n\right)\geq 1 -
n^{-K},\quad\forall\, K>0, ~n\geq n(K).
\]
\end{lemma}
\noindent{\bf  Remark.~} Borrowing a term from Knuth et al. \cite{KMP},
the event $\Delta$ holds quite surely (q.s.).

\begin{proof} Introduce the events
\[
C_t=\bigg\{B(s+1)\leq
\Psi_\mbp(B(s))~\forall~s<t,~B(t+1)>\Psi_\mbp(B(t)), B(t) \geq
k_*\bigg\};
\]
that is, $C_t$ is the event that the recursive inequality $B(s+1)\le
\Psi_{\mbp}(B(s))$ is violated at a state $B(t)\ge k_*$, and $t$ is
the first such moment. Clearly
\[
\Delta^c=\bigcup_{t\ge 0}C_t.
\]
Let us show that
\[
C_t=\emptyset,\quad\forall\,t\ge t_*:= 4c_2^{-1/2}\ln n.
\]
Suppose that on the contrary $C_{t_1}\neq\emptyset$ for some
$t_1>t_*$. Then, by the definition of $C_{t_1}$, we have
\[
B(s+1)\le \Psi_{\mbp}(B(s)),\quad\forall\, s\le t_*,
\]
and certainly $B(t_*)\ge k_*$. However, using this recurrence
inequality exactly as in the derivation of \eqref{solution}, we must
have
\begin{align*}
B(t_*) &\leq n\left(1-c_2^{1/2}/4\right)^{t_*} +2c_2^{-1/2}\\
&\leq n\exp\left[-\left(4c_2^{-1/2}\right)(\ln n)  \left(c_2^{1/2}/4\right)\right] +2c_2^{-1/2}\\
   &= 1+2c_2^{-1/2}
<k_*
\end{align*}
since $c_2=o\left( \ln^{-2\ep} n\right)$. Contradiction!  Thus $\Delta^c$ is a
union of at most $t_*$ events $C_t$. Now by \eqref{e100} and Lemma
\ref{decr}, we have
\[
P(C_t)\le 3n^{1/2}\exp(-H_{\mbp}(k_*)/2),\quad t\ge 0,
\]
since $B(t)\ge k_*$ on $C_t$. Therefore
\begin{align*}
P(\Delta^c)=P\bigg(\bigcup_{t\le t_*}C_t\bigg)&\le \sum_{t\le
  t_*}P(C_t)\\ &\le
3t_*n^{1/2}\exp\big(-H_{\mbp}(k_*)/2\big)\\ 
&=12c_2^{-1/2}n^{1/2}\ln n\exp\big(-H_{\mbp}(k_*)/2\big)\\ 
&\le 12n\ln n\exp\big(-H_{\mbp}(k_*)/2\big)\\ 
&\le 12n\ln n\exp\big(-0.5A_*\ln^{1+\varepsilon}n\big) ~~~~~~~\mbox{by Lemma \ref{Hlemma}},
\end{align*}
and from here the lemma follows.
\end{proof}
This completes a program we put forth at the end  of the introduction. Combining Lemma \ref{deltalemma} and this last Lemma \ref{lemma3.8}, we have proved the following.

\begin{lemma}\label{lemma3.9} Let $k_*=c_2^{-1}\ln^{-\ep} n$ for some $\ep \in (0,1/4)$ however small, and let $\tau(k_*)$ denote the random moment when $B(t)$ falls to or below $k_*$ for the first time. Then, for some constant $a>0$,
\[
P\bigl\{\tau(k_*)\le ac_2^{-1/2}\ln n\bigr\}\geq 1 -
n^{-K},\quad\forall K>0,~n\geq n(K).
\]
\end{lemma}
\noindent In short, q.s. $\tau(k_*)=o(c_2^{-1})$.

\subsection{Bounding the expectation of $\tau(k_*)$}
Even though $\tau(k_*)=o(c_2^{-1})$ q.s., proving that
$E[\tau(k_*)]=o(c_2^{-1})$ as well is not straightforward, since we do
not have a polynomial (worst-case) bound for $\tau(k_*)$.  As a first
step, introducing the event indicators $I_{\Delta}$, $I_{\Delta^c}$,
we split $E[\tau(k_*)]$ using $1=I_{\Delta}+I_{\Delta^c}$ and bound
the second summand via the Cauchy-Schwarz inequality:
\begin{eqnarray}\notag
E[\tau(k_*)] =& E[\tau(k_*)I_{\Delta}]&+~E[\tau(k_*)I_{\Delta^c}
]\\ \notag
\leq & 5(\ln n)c_2^{-1/2} &+ ~\sqrt{ E[\tau^2(k_*)] \cdot
  E[I_{\Delta^c}^2]}\\ \label{CS} \leq & 5(\ln n)c_2^{-1/2} &+ ~\sqrt{
  E[\tau^2(k_*)]}\cdot \sqrt{ P(\Delta^c)}.
\end{eqnarray}
It remains to show that $E[\tau^2(k_*)]$ is at most polynomially
large. To do so, introduce $T_k$, the random time the process
$\{B(t)\}$ spends at state $k$, i.e.
\begin{equation}\label{e68}
T_k=|\{t\ge 0\,:\,B(t)=k\}|.
\end{equation}
Then
$$ 
\tau^2(k_*)=\left(\sum_{k=k_*+1}^n T_k\right)^2\le
(n-k_*)\sum_{k=k_*+1}^n T_k^2,
$$
again by the Cauchy-Schwarz inequality. Therefore
$$
 E[\tau^2(k_*)]\,\le\, n\!\!\!\!\sum_{k=k_*+1}^n E[T_k^2]\,\,\le
n\!\!\!  \sum_{k=k_*+1}^n E[T_k^2\,|\,T_k>0].
$$
Recalling the notation $\pi_{kk}=P(B(t+1)=k\,|\,B(t)=k)$, we
observe that
$$
 P\{T_k=j\,|\,T_k>0\} = \pi^{j-1}_{kk} (1-\pi_{kk}),\quad j>0,
$$
 i.e. conditioned on $\{T_k>0\}$, $T_k$ is geometrically
distributed, with success probability $1-\pi_{kk}$. In particular,
$$
 E[T_k\,|\,T_k>0]=\frac{1}{1-\pi_{kk}},\quad
\text{Var}[T_k\,|\,T_k>0]=\frac{\pi_{kk}} {(1-\pi_{kk})^2}.
$$
It is obvious intuitively, and can be easily proved, that
$\pi_{kk}$ decreases with $k$: the larger the number of balls---the
larger the probability of collision. Consequently both conditional
moments of $T_k$ decrease with $k$. So
\begin{align}\notag
E[\tau^2(k_*)] &\leq n\sum_{k=k_*+1}^n
\left(\frac{\pi_{kk}}{(1-\pi_{kk})^2} +
\frac{1}{(1-\pi_{kk})^2}\right)\\ \notag &\leq n\sum_{k=k_*+1}^n
\frac{2}{(1-\pi_{kk})^2} \\ &\leq
\frac{2n^2}{(1-\pi_{k_*k_*})^2}. \label{etau}
\end{align}
Therefore, it remains to show only that $\frac{1}{1-\pi_{k_*k_*}}$ is at most polynomially  large in $n$.
Using the simplest {\it lower\/} bound for the probability of the union of events,  via the inclusion-exclusion formula, we write
\begin{align}\notag
1-\pi_{kk} &= P(\mbox{there is a collision during a   $k$-allocation})\\
\notag &\geq \sum_{\{a,b\}\subset [k]} P(\mbox{balls $a$ and $b$  collide})-\!\!\!\!\!\!\!\sum_{\{c,d\}\neq\{e,f\}\subset[k]}P(\mbox{$c,d$  collide and $e,f$ collide})\\ \label{ie0}
 & \geq \binom{k}{2} c_2 - \left(\binom{k}{3}c_3 + \frac{1}{2}\binom{k}{2}\binom{k-2}{2} c_2^2\right)\\\label{inclexcl}
&=\binom{k}{2} c_2\left(1-a_1k\frac{c_3}{c_2}-a_2k^2 c_2\right).
\end{align}
where $a_1, a_2$ are some absolute constants. The two rightmost terms of \eqref{ie0} are due to the fact that there are two ways in which two distinct pairs of balls can collide: either the two pairs overlap at one ball or they are disjoint. Now introduce
$k_1:=c_2^{-1/2}\ln^{-\ep/4} n$; clearly $k_1<k_*$, and we also have
$$ 
\lim k_1\frac{c_3}{c_2}= 0,\quad \lim k_1^2c_2=0.
$$
So, by \eqref{inclexcl}, uniformly for all $k\le k_1$,
\begin{equation}\label{e71}
\frac{1}{1-\pi_{kk}}\le \frac{1}{\binom{k}{2}c_2}(1+o(1)).
\end{equation}
Consequently
\begin{equation}\label{e72}
\frac{1}{1-\pi_{k_*k_*}}\le
\frac{1}{1-\pi_{k_1k_1}}=O(k_1^{-2}c_2^{-1})= O(\ln^{\ep/2} n),
\end{equation}
which is sufficient. Thus by \eqref{CS} and \eqref{etau},
\begin{equation}\label{e73}
E[\tau(k_*)]=O(c_2^{-1/2}\ln n)=o(c_2^{-1}).
\end{equation}

\section{Bounding expected duration of a middle phase and a late phase}\label{midlate}

To complete the proof of Theorem \ref{t1} it remains to bound the
expected duration of the process after the number of balls has dropped
below $k_*$.

We define a middle phase as $[\tau(k_*),\tau(k_1))$, the time interval
  during which the number of balls is below $k_*$ and above
  $k_1$. Using the $T_k$ defined in \eqref{e68}, we have
\begin{equation}\notag
E[\tau(k_1)-\tau(k_*)]=\sum_{k=k_1+1}^{k_*}E[T_k]\le
\sum_{k=k_1+1}^{k_*}\frac{1}{1-\pi_{kk}}.
\end{equation}
Then, by decreasing monotonicity of $(1-\pi_{kk})^{-1}$
and \eqref{e72},
\begin{equation}
E[\tau(k_1)-\tau(k_*)]\le
\frac{k_*-k_1}{1-\pi_{k_1k_1}}=O\bigl(c_2^{-1}\ln^{-\ep+\ep/2}n\bigr)=
O\bigl(c_2^{-1}\ln^{-\ep/2}n\bigr)=o(c_2^{-1}).\label{e201}
\end{equation}
(The last computation explains at long last why we needed the
$\ln^{-\ep} n$ factor in the definition \eqref{kstardef} of $k_*$.)

Naturally, we define a late phase as $[\tau(k_1),\tau(1)]$. By
\eqref{e71},
\begin{equation}
E[\tau(1)-\tau(k_1)]\le
1+\sum_{k=2}^{k_1}\frac{1+o(1)}{c_2\binom{k}{2}}=2c_2^{-1}(1+o(1)),\label{e75}
\end{equation}
where we used
$$ \frac{1}{k(k-1)}=\frac{1}{k-1}-\frac{1}{k},\quad k\ge 2.
$$

That does it! Adding the estimates \eqref{e73}, \eqref{e201},
and\eqref{e75}, we obtain
$$ 
E[\tau(1)]\le 2c_2^{-1}(1+o(1)).
$$

\section{Lower Bound}\label{lb}
We  now  provide the matching lower bound for Theorem \ref{t1}. We will not need assumptions on  $\mbp$  as strong  as \eqref{conds}; rather we will simply  assume that $c_2(\mbp)\to 0$. 

For any $m\leq n$, let $T(m)$ denote the coalescence time for the process starting with $m$ balls (we have so far been considering the case $T:=T(n)$.) We start by stating Theorem 2  from \cite{AAKR}:
 \begin{equation}\label{AAKRt2}
\mbox{For any } \mbp,  ~~~~ E[T(m)] \geq 2c_2^{-1}\left(1-\frac{1}{m} -\frac{(m-1)(m-2)}{12}\frac{c_3}{c_2}\right).
\end{equation}
To get a bound on $T(n)$, first note the obvious-looking fact that  
\begin{proposition}\label{coupling}
Let $m_1 \leq m_2 \leq n$. Then for any $\mbp$, $T(m_2)$ stochastically dominates $T(m_1)$. 
\end{proposition}
\begin{proof}
This is a result of the  following basic coupling  argument: start with $m_2$  balls, $m_1$ of which are marked. Then perform the usual allocation  process; the time $X$ at which all balls coalesce is distributed as $T(m_2)$, and at this time certainly all marked balls have coalesced  as well; call the time that these marked balls have first coalesced  $Y$,  so that $X\geq Y$ and $Y$ is distributed  as $T(m_1)$. From here the result follows.
\end{proof}
This result shows that for the bound  $E[T(n)]\geq 2c_2^{-1}(1-o(1))$, it is sufficient to show $E[T(m)]\geq 2c_2^{-1}(1-o(1))$  for some $m\leq n$. By the inequalities
\[
 c_2^2\leq c_3 \leq c_2^{3/2}~~\Longrightarrow~~c_2\leq \frac{c_3}{c_2} \leq c_2^{1/2},
 \]
we have that $c_2 \to 0$  iff  $c_3/c_2\to 0$. Let $m_*:=(c_2/c_3)^{1/3}$, say; then $m_* \to \infty$ since we assume that $c_2\to 0$, and so \eqref{AAKRt2} becomes
\[
  E[T(m_*)] \geq 2c_2^{-1}\left(1-\frac{1}{m_*} -O(m_*^2m_*^{-3})\right)=2c_2^{-1}(1-o(1)),
\]
and the result  follows by Proposition \ref{coupling}. The proof Theorem \ref{t1} is complete.

\section{Proof of Theorem \ref{cexthm}}\label{cexsec}

We start by restating Theorem \ref{cexthm} in a more detailed manner.
\begin{theorem}\label{t5.1} 
Let $b_0=n$. Suppose $c_2=(\ln^{-2} n)\lambda(n)$, where $\lambda(n) \to \infty$ however slowly. Set
$\mbp=\mbta(c_2)$, i.e. $\mbp$ is the topheavy
distribution $(\theta_1,\theta_2,\dots,\theta_2)$ with
$$
\theta_1^2 + (n-1)\theta_2^2=c_2.
$$
 Then for the process evolving according to this $\mbp$,  whp
$$
\tau(1) \ge c_2^{-1}\frac{\sqrt{\lambda(n)}}{20};
$$
so, in particular, $E[\tau(1)]=\omega\left( c_2^{-1}\right)$.
\end{theorem}
\noindent{\bf Remarks.}
\begin{enumerate}
\item Our choice of $\mbp$ should be expected. Indeed, the recurrence
  inequality \eqref{thetacor1} signals, intuitively, that the
  coalescent process for $\mbp=\mbta(c_2)$ is a good candidate for
  being the slowest among all $\mbp$ with $\sum_jp_j^2 =c_2$.

\item In Theorem 5 of \cite{AAKR}, Adler et al. had proved that for
  $\max_j p_j$ bounded away from $0$, and the remaining $p_j$
  uniformly small, the expected coalescence time exceeds $c_2^{-1}$ by
  a factor of $\ln n$.
\end{enumerate}

\begin{proof}[Proof of Theorem \ref{t5.1}]
For simplicity we let $c:=c_2$.  For $\mbp=\mbta(c)$, the
Markov chain is almost as simple as that for the uniform
$\mbp$. Indeed, given $B(t)$, the number of balls that land in box
$2,\dots, n$ (call it $\hat B(t)$) is binomially distributed with
parameters $B(t)$ and success probability $1-\theta_1$, i. e. $\hat
B(t)=\text{Bin}(B(t),1-\theta_1)$ in short.  Conditioned on $\hat
B(t)$, we have a uniform allocation of $\hat B(t)$ balls among $n-1$
boxes $2,\dots,n$. And, for $B(t)$ sufficiently large,  whp $\hat
B(t)\sim (1-\theta_1)B(t)$.  So, based on our experience with
deterministic approximations earlier in the paper, we should expect
that---after fusing balls that landed in the same box---these $\hat
B(t)$ balls give birth to about
$$ (n-1)\left[1-\exp\left(-\frac{\hat B(t)}{n-1}\right)\right]\sim
(n-1)\left[1-\exp\left(-\frac{(1-\theta_1) B(t)}{n-1}\right)\right]
$$ balls for next generation. All the balls that landed in box $1$, if
there are any, will coalesce into one ball. Ignoring this box for now,
we expect then that the process $\{B(t)\}$ whp ``closely'' obeys a
recurrence inequality of the form
\begin{equation}\notag
B(t+1)\ge (n-1)\varphi\left(\frac{(1-\theta_1) B(t)}{n-1}\right),\quad
\varphi(x):= 1-e^{-x}.
\end{equation}
Here is a precise claim.
\begin{lemma}\label{l5.2}
 Let $\gamma := 1-2c^{1/2}<1-\ta_1$ and introduce
$$ \Psi(k)=(n-1)~\eta\left(\frac{\gamma k}{n-1}\right),\quad
 \eta(x):=1.5\varphi(x)-0.5 x.
$$ Then, for $n$ sufficiently large,
\begin{equation}\label{e75b}
P(B(t+1) < \Psi(B(t))\,|\,B(t)=k)\le e^{-kc/3}+e^{-k^3/73n^2}.
\end{equation}
\end{lemma}
\noindent{\bf Remark.~} $\eta(0)=0$, $\eta(x) \leq \varphi(x)$, and $\eta(x)$ is increasing for
$x\le \ln 3$.

\begin{proof}[Proof of Lemma \ref{l5.2}] Notice first that 
$$ B(t+1)\ge \mathcal{R}(t+1),
$$ where $\mathcal{R}(t+1)$ is the number of boxes among $2,\dots,n$
  that host at least one of $\hat B(t)$ balls.  Denoting
  $P(\{\cdot\}|B(t)=k)$ by $P_k(\{\cdot\})$, we have then
\begin{align*}
P_k(B(t+1)<\Psi(k))\le& P_k(\mathcal{R}(t+1)<\Psi(k))\\ \le&
P_k(\mathcal{R}(t+1)<\Psi(k),\, \hat B(t)\ge \gamma k) +P_k(\hat
B(t)<\gamma k).
\end{align*}
Now, denoting the c.d.f. of $\mathcal{R}(t+1)$ conditioned on $\{\hat
B(t)=j\}$ by $F_j$, we have: for $j_1<j_2$,
$$ F_{j_2}(x) \le F_{j_1}(x),\quad\forall x\ge 0.
$$ (Informally, the fewer balls we allocate among the boxes
$2,\dots,n$, the fewer nonempty boxes we end up with.) Therefore
\begin{align*}
P_k(\mathcal{R}(t+1)<\Psi(k),\, \hat B(t)\ge \gamma k) =&\sum_{j\ge
    \gamma k} P(\mathcal{R}(t+1)<\Psi(k)\,|\,\hat B(t)=j)\,P_k(\hat
  B(t)=j)\\ \le&P(\mathcal{R}(t+1)<\Psi(k)\,|\,\hat B(t)=\lceil \gamma
  k\rceil)\sum_{j\ge \gamma k} P_k(\hat B(t)=j)\\ \le&
  P(\mathcal{R}(t+1)<\Psi(k)\,|\,\hat B(t)=\lceil \gamma k\rceil).
\end{align*}
Consequently
\begin{align*}
P_k(B(t+1)<\Psi(k)) \le& \mathcal{P}_1
+\mathcal{P}_2,\\ \mathcal{P}_1:=&
P_k(\mathcal{R}(t+1)<\Psi(k)\,|\,\hat B(t)=\lceil\gamma
k\rceil),\\ \mathcal{P}_2:=&P_k(\hat B(t)<\gamma k).
\end{align*}
Using the Chernoff bound for the tail of binomial distribution (see
Mitzenmacher and Upfal \cite{MU}, for instance), we have
\begin{align}\notag
\mathcal{P}_2=&P(\text{Bin}(k,1-\theta_1)<\gamma k)\\\notag
\le&\left.\exp\left(-\frac{k p\delta^2}{2}\right)\right|_{
  p=1-\theta_1,~ \delta= 1-\gamma
  (1-\theta_1)^{-1}}\\ \le&e^{-kc/3},\label{e76}
\end{align}
as $\theta_1\sim c_2^{1/2}$, from the definition
\eqref{thetaproperties} of $\ta_1,\ta_2$. 

Turn to $\mathcal{P}_1$. Applying Theorem \ref{Chern1} to the boxes set
$\{2,\dots,n\}$ and the uniform distribution $\bold u$ on this set, we
obtain
$$ \mathcal{P}_1\le 3\sqrt{\gamma k}\exp\left(-\frac{(\Phi_{\bold
    u}(\gamma k)-\Psi(k))^2}{2\gamma k}\right).
$$ Here, using the definition of $\eta(x)$,
\begin{align*}
\Phi_{\bold u}(\gamma k)-\Psi(k)=&(n-1)\varphi\left(\frac{\gamma
  k}{n-1}\right)-(n-1) \eta\left(\frac{\gamma
  k}{n-1}\right)\\ =&\frac{(n-1)}{2}\left[\frac{\gamma
    k}{n-1}-\varphi\left(\frac{\gamma k}{n-1}\right)\right].
\end{align*}
So, as
$$ x-\varphi(x)=e^{-x}-1+x \ge \frac{x^2}{2}-\frac{x^3}{6}\ge
\frac{x^2}{3},\quad x\in [0,1],
$$ we have
\begin{equation}\label{e77}
\mathcal{P}_1\le 4\sqrt{\gamma k}\exp\left(-\frac{(\gamma
  k)^4}{72\gamma k (n-1)^2}\right)\le e^{-k^3/73n^2}.
\end{equation}
The estimates \eqref{e76}-\eqref{e77} imply \eqref{e75b}, thereby completing the proof of the lemma.
\end{proof}

To continue, let $\hat k:= n^{3/4}$. Introduce two events,
$$ \Gamma:=\{\forall t,\,B(t)\ge \hat k \Longrightarrow B(t+1)\ge
\Psi(B(t))\},
$$ and
$$ \Pi:=\{\exists t > c^{-3/2}\,:\,B(t)\ge \hat k\}.
$$ On the event $\Gamma$, $B(t)$ does not decrease ``too quickly'' as
long as $B(t)$ is above $\hat k$. On the event $\Pi$, $\tau(\hat k)\ge
c^{-3/2}$, i. e. $\tau(\hat k)=\omega\left(c^{-1}\right)$.  Then, by Lemma \ref{l5.2},
\begin{align*}
P(\Gamma^c\cap\Pi^c)\le&P\left\{\bigcup_{t=0}^{c^{-3/2}}\{B(t+1)<\Psi(B(t)),\,B(t)
\ge \hat k\}\right\}\\ \le&(1+c^{-3/2})\left(e^{-\hat k c/3}+e^{-\hat
  k^3/73
  n^2}\right)\\ =&(1+c^{-3/2})\left(e^{-n^{3/4}c/3}+e^{-n^{1/4}/73}\right)\to
0,
\end{align*}
i. e. $P(\Gamma\cup \Pi)\to 1$.

If we show that $\tau(\hat k)=\omega\left( c^{-1}\right)$ on the event $\Gamma$ as
well, we will be able then to claim that  whp $\tau(\hat k)=\omega\left( c^{-1}\right)$,
and the proof of Theorem \ref{t5.1} will be complete. 

To do so, we observe that on the event $\Gamma$,
\begin{equation}\label{e79}
x(t+1)\ge \eta(\gamma x(t)),\quad x(t):=\frac{B(t)}{n-1},\quad
x(0)=1-(n-1)^{-1},
\end{equation}
as long as
$$ x(t)\ge \frac{\hat k}{n-1}\sim n^{-1/4}.
$$
\begin{lemma}\label{l5.3} Let $n\ge 3$. Under the recurrence \eqref{e79},
\begin{equation}\label{e80}
x(t)\ge \frac{2}{3}\,\frac{\gamma^t}{t+1}.
\end{equation}
\end{lemma}
\begin{proof} The argument runs in parallel  to that for the {\it lower\/}  bound of $B(t)$
in \eqref{eq1}. The base case $t=0$ is just
$$ x(0)=1-n^{-1}\ge \frac{2}{3}.
$$ Suppose \eqref{e80} holds for some $t$. Since $\eta(x)$ is
increasing for $x\le \ln 3$, \eqref{e79} implies that
$$ x(t+1)\ge\eta\left(\gamma\frac{2}{3}\frac{\gamma^t}{t+1}\right)=
\eta\left(\frac{2}{3}\frac{\gamma^{t+1}}{t+1}\right).
$$ So, to complete the inductive step, we need to show that
\begin{equation}\label{e81}
\eta(y)\ge y\,\frac{t+1}{t+2},\quad
y:=\frac{2}{3}\frac{\gamma^{t+1}}{t+1}.
\end{equation}
Define $z$ as a root of
\begin{equation}\label{e82}
\eta(z)=z\frac{t+1}{t+2}\quad\text{or}\quad 1-e^{-z}
=z\frac{t+4/3}{t+2}.
\end{equation}
Since $1-e^{-1}< 2/3$, equation \eqref{e82} has a (unique) root
$z=z(t)$ for $t\ge 0$. Using $1-e^{-z}\ge z-z^2/2$, we obtain
$$ z(t)\ge \frac{4/3}{t+2}.
$$ Inequality \eqref{e81} holds if $y\le z(t)$, which is certainly so
because
$$
y(t) \leq \frac{2/3}{t+1}\le \frac{4/3}{t+2},\quad \forall\,t\ge 0,
$$
and  from  here the lemma follows.
\end{proof}

Thus on the event $\Gamma$,
\begin{equation}\label{e83}
B(t)\ge \hat k\Longrightarrow B(t+1)\ge
(n-1)\frac{2}{3}\frac{\gamma^{t+1}}{t+2}.
\end{equation}

\begin{lemma}\label{l5.4} On the event $\Gamma$,
\begin{equation}\notag
\tau(\hat k) \ge \frac{\ln^2 n}{20\sqrt{\lambda(n)}}=c^{-1}\frac{\sqrt{\lambda(n)}}{20}.
\end{equation}
\end{lemma}

\begin{proof}[Proof of Lemma \ref{l5.4}] Let the event $\Gamma$ hold. By \eqref{e83}, and the definition of
$\tau(\cdot)$, $\hat\tau:=\tau(\hat k)$ satisfies
$$ \hat k\ge B(\hat\tau)\ge
  (n-1)\frac{2}{3}\frac{\gamma^{\hat\tau}}{\hat\tau+1}.
$$ Recalling that $\hat k=n^{3/4}$, $\gamma = 1-2c^{1/2}$, and taking
  logarithms, we get
$$ 3c^{1/2}\hat\tau+\ln(\hat\tau+1)\ge \frac{1}{5}\ln n,
$$ or
\begin{equation}\label{tauhat}
3\frac{\sqrt{\lambda(n)}}{\ln n}\hat\tau+ \ln(\hat\tau+1)\ge
  \frac{1}{5}\ln n.
\end{equation}
 It then follows that for large enough $n$, 
$$
\hat\tau\ge \frac{1}{4\cdot 5}\frac{\ln^2 n}{\sqrt{\lambda(n)}}=c^{-1}\frac{\sqrt{\lambda(n)}}{20},
$$
as if it did not, then $\ln(\hat \tau +1) =o(\ln n)$ and 
\[
3\frac{\sqrt{\lambda(n)}}{\ln n}\hat\tau \leq \frac{3}{4\cdot 5}\ln n <\ln n,
\]
which together contradict \eqref{tauhat} for $n$ large enough. This completes the proof of Lemma \ref{l5.4}.
\end{proof}

In summary, on the event $\Gamma\cup \Pi$,
$$ \tau(\hat k)\ge \min\left\{c^{-3/2},c^{-1}\frac{\sqrt{\lambda(n)}}{20}\right\}=c^{-1}\frac{\sqrt{\lambda(n)}}{20}.
$$ Recalling that $P(\Gamma\cup\Pi)\to 1$, we conclude that
\begin{align*}
c E[\tau(1)]\ge c E[\tau(\hat k)]\ge E[\tau(\hat k)\,I_{\Gamma\cup
    \Pi}] \ge \frac{\sqrt{\lambda(n)}}{20}\,P(\Gamma\cup\Pi)\to\infty.
\end{align*}
This concludes the proof of Theorem \ref{t5.1}.
\end{proof}

\section{Concluding remarks, future work}\label{fut}
  
A generalization of this problem is to allow the probability of a ball
going to a certain box to depend on its \textit{origin} and not just on
its \textit{destination}; that is, if we have a ball in box $i$, it
has probability $p_{ij}$ of landing in box $j$ for any $j \in [n]$,
and these probabilities are not necessarily the same for all $i$. This
is more difficult, as $\{B(t)\}_{t\geq 0}$ is no longer a Markov
chain: we have to keep track of the \textit{locations} of the balls at
any time and not simply their number.  

Coupling from the past algorithms involve running simultaneous coalescing flows on a Markov chain  $P:=(p_{ij})_{1\leq i,j \leq n}$ with stationary distribution $\mbpi$, and return samples distributed {\it exactly} according to $\mbpi$, at the time when all the  flows coalesce. It would be very interesting to extend the techniques in this paper to a more general case (when the rows of $P$ are not necessarily all equal to some vector $\mbp$) in order to obtain an {\it upper} bound for the expected running time of such algorithms. \\

Here is an approach that appears promising. We start with a fixed allotment of one ball in each of $n$ boxes, and for $t_*$ times we run $n$ independent allocations for each ball, where we do not fuse balls that land in the same box; call this a {\it mixing phase}. This $t_*$ is to be taken large enough so that the location of each ball at time $t_*$ is ``almost'' $\mbpi$-distributed. Then at time $t_*+1$, we  allocate the balls into the boxes but fuse any that collide (a {\it fusing allocation}). Continue alternating between mixing phases and fusing allocations until total coalescence has occurred; the coalescence time for this process should dominate the time for the usual process where we fuse colliding balls at every time. Moreover, if the locations of the balls at each fusing allocation are sufficiently independent and close to $\mbpi$-distributed, then we may be able to bring the results in this paper to bear on this more general case.


\appendix
\section{Missing parts of the proof of Theorem \ref{Chern1}}
\begin{lemma}\label{hess}
Let $\chi= H_{rr}H_{zz}-H_{rz}^2$. Then $\chi>0$ for all $(z,r,b)$ on
the curve $C$ defined by \eqref{hz} and \eqref{hr}.
\end{lemma}
\begin{proof}
By the definition of $H$, we have
\begin{equation}\label{hzz}
H_{zz} = \frac{b}{z^2} - \sum_j
\frac{(e^{np_jr}-1)^2}{(1+z(e^{np_jr}-1))^2},
\end{equation}
\begin{equation}\label{hrr}
    H_{rr}= \frac{k}{r^2} +z(1-z)n^2 \sum_j \frac{p_j^2
      e^{np_jr}}{(1+z(e^{np_jr}-1))^2},
\end{equation}
and
\begin{equation}\label{hrz}
H_{rz}=n\sum_j \frac{p_j
  e^{np_jr}}{(1+z(e^{np_jr}-1))^2}.
\end{equation}
We can recast \eqref{hz} and \eqref{hr} as
\begin{align}\label{hz7}
b = \sum_j \frac{z(e^{np_jr}-1)}{1+z(\enp-1)},~~~~ k =nrz\sum_j
\frac{p_j \enp }{1+z(\enp-1)}.
\end{align}
Using these in \eqref{hzz}-\eqref{hrz}, we get that on $C$,
\begin{equation}\label{hzz3}
H_{zz} = \frac{1}{z^2}\sum_j \frac{z(e^{np_jr}-1)}{1+z(\enp-1)}
-\sum_j\frac{(e^{np_jr}-1)^2}{(1+z(e^{np_jr}-1))^2},
\end{equation}
\begin{equation}\notag
H_{rr} = \frac{nz}{r}\sum_j \frac{p_j \enp }{1+z(\enp-1)} +z(1-z)n^2
\sum_j \frac{p_j^2\enp}{(1+z(\enp-1))^2},
\end{equation}
and
now \eqref{hzz3} simplifies to
\begin{equation}\label{hzz4}
H_{zz}= \frac{1}{z} \sum_j \frac{(\enp -1)}{(1+z(\enp -1))^2}.
\end{equation}
For $H_{rr}$ we can put under a common denominator and get
\begin{equation}\notag
    H_{rr} = \sum_j
    \frac{nzp_j\enp(1+z(\enp-1))+rz(1-z)n^2p_j^2\enp}{r(1+z(\enp-1))^2}.
\end{equation}
Using the inequality $e^x-1\geq x$ gives
\begin{equation}\notag
H_{rr} \geq \sum_j
\frac{nzp_j\enp(1+znp_jr)+rz(1-z)n^2p_j^2\enp}{r(1+z(\enp-1))^2},
\end{equation}
which then leads to some very convenient cancelling (in particular, of
the $z^2$ term) to get
\begin{equation}\label{hrr6}
    H_{rr} \geq \sum_j \frac{nzp_j\enp (1+rnp_j)}{r(1+z(\enp-1))^2}.
\end{equation}
Multiply \eqref{hzz4} and \eqref{hrr6} together:
\begin{equation}\notag
    H_{zz}H_{rr} \geq \frac{n}{r} \left( \sum_j \frac{(\enp
      -1)}{(1+z(\enp -1))^2} \right) \left(\sum_j \frac{p_j\enp
      (1+rnp_j)}{(1+z(\enp-1))^2} \right),
\end{equation}
and use the Cauchy-Schwarz inequality to get
\begin{equation}\label{hzzhrr2}
H_{zz}H_{rr} \geq \frac{n}{r} \left( \sum_j
\frac{[(\enp-1)p_j(1+np_jr)\enp]^{1/2}}{(1+z(\enp-1))^2} \right)^2.
\end{equation}
We need to show that this is strictly greater than $H_{rz}^2$, which
can be expressed, using \eqref{hrz}, as
\begin{equation}\label{hrzs}
H_{rz}^2 = n^2 \left(\sum_j \frac{p_je^{np_jr}}{(1+z(e^{np_jr}-1))^2}
\right)^2.
\end{equation}
Taking square roots of the expressions \eqref{hzzhrr2} and
\eqref{hrzs}, the condition $\chi>0$ is equivalent to
\begin{equation}\label{e100a}
\sum_j \frac{[(\enp-1)p_j(1+np_jr)\enp]^{1/2}-\sqrt{rn}p_j
  \enp}{(1+z(e^{np_jr}-1))^2} > 0.
\end{equation}
For \eqref{e100a} to hold, it suffices that each summand is
nonnegative (and at least one strictly positive). Multiplying the
numerators by $\sqrt{rn}$, we need to show that $\forall~j,~1\leq j
\leq n$,
\begin{equation}\notag
[(\enp-1)np_jr(1+np_jr)\enp]^{1/2}-rnp_j \enp \geq 0.
\end{equation}
This is equivalent to showing that $\forall~j$,
\begin{equation}\label{npr}
f(np_jr)\geq 0,
\end{equation}
where we define $f(x)=(e^x-1)x(1+x)e^x - x^2e^{2x}$. Now
\begin{align*}
f(x)&= xe^{2x} +x^2e^{2x} -xe^x -x^2e^x -x^2e^{2x}\\ &= xe^{2x}-xe^x
-x^2e^x= xe^x(e^x-1-x) > 0~~~~~~\mbox{for}~~x>0.
\end{align*}
Therefore the inequalities in \eqref{npr} hold (and at least one of
them is strict), and so the lemma follows.
\end{proof}
\begin{lemma}\label{hdplemma}
With $h(b)$ as defined in \eqref{hdef}, we have, uniformly for  $b\leq k$,
\[
h''(b) \leq -\frac{1}{k}.
\]
\end{lemma}
\begin{proof}
First note that 
\begin{align}\notag
h''(b) =\frac{d}{db} \left(h'(b)\right) &=\frac{d}{db}
\left(-\ln(z(b))\right)~~~~~~~~~\mbox{by \eqref{hprime}}\\
&=-z'(b)/z(b).\label{z'z}
\end{align}
To find $z'(b)$, differentiate $H_z(z(b),r(b),b)=0$ and $H_r(z(b),r(b),b)=0$ with respect to $b$; we can solve  for $z'(b)$  in this system (using $H_{rb}=0, H_{zb}=-1/z$) to get
\begin{equation}\label{derz}
z'(b)= \frac{H_{rr}}{z(b) (H_{zz}H_{rr}-H_{rz}^2)}= \frac{H_{rr}}{z(b)\chi}, 
\end{equation}
which is strictly  positive  by \eqref{hrr6} and Lemma \ref{hess}. Now note that, using the expression \eqref{hzz4} for $H_{zz}$, and \eqref{hz7} to express $b$, 
\begin{align}\notag
\frac{H_{zz}}{b}  &= \frac{1}{z(b)^2b} \sum_j \frac{z(b)(e^{np_jr(b)}-1)}{[1+z(b)(e^{np_jr(b)}-1)]^2}\\
& \leq \frac{1}{z(b)^2b} \sum_j \frac{z(b)(e^{np_jr(b)}-1)}{1+z(b)(e^{np_jr(b)}-1)} =\frac{1}{z(b)^2}. \label{e110}
\end{align}
Therefore, using \eqref{z'z}, \eqref{derz}, and \eqref{e110},
\begin{align*}
h''(b) =-\frac{z'(b)}{z(b)}  = -\frac{H_{rr}}{z(b)^2\chi} \leq \frac{-H_{rr}H_{zz}}{b\chi} \leq -\frac{1}{b},
\end{align*}
where  the last inequality holds because $0<\chi <H_{rr}H_{zz}$. Using the fact that $b \leq k$, $\Phi_\mbp(k) \leq k$, and $\tilde{b}$ is between $b$ and $\Phi_\mbp(k)$, we obtain
\[
h''(\tilde{b}) \leq -\frac{1}{k}.
\]
\end{proof}

\section*{Acknowledgments}
We are very grateful to the participants of a student workshop on combinatorial probability at Ohio State University
for many productive discussions of this project. We would  also like to thank an anonymous reviewer for many insightful corrections and suggestions that greatly helped us improve the presentation of the paper.


%
%
%
%

\end{document}